# A STABILIZED MIXED FORMULATION FOR *UNSTEADY* BRINKMAN EQUATION BASED ON THE METHOD OF HORIZONTAL LINES

SHRIRAM SRINIVASAN AND K. B. NAKSHATRALA

ABSTRACT. In this paper, we present a stabilized mixed formulation for *unsteady* Brinkman equation. The formulation is systematically derived based on the variational multiscale formalism and the method of horizontal lines. The derivation does not need the assumption that the fine-scale variables do not depend on the time, which is the case with the conventional derivation of multiscale stabilized formulations for transient mixed problems. An expression for the stabilization parameter is obtained in terms of a bubble function, and appropriate bubble functions for various finite elements are also presented. Under the proposed formulation, equal-order interpolation for the velocity and pressure (which is computationally the most convenient) is stable. Representative numerical results are presented to illustrate the performance of the proposed formulation. Spatial and temporal convergence studies are also performed, and the proposed formulation performed well.

## 1. INTRODUCTION

Among the engineering and environmental issues that demand attention at present, enhanced oil recovery (EOR), carbon-dioxide sequestration, and contaminant transport in heterogeneous media can lay claim to pride of place. The models currently used are mostly steady state models and reliable computational techniques do exist for them. However, there is a big drawback because it is only if the transient effects are considered that it will be possible to later study the coupling phenomena that occur due to interaction between different media. Our intention is to propose a stabilized mixed formulation for an unsteady model, since the existing techniques have their own disadvantages. However, before developing the finite element formulation, one needs to turn to theory of mixtures to trace the development of the models.

The foundations of the theory of interacting continua (also known as theory of mixtures) were laid down by Truesdell [1, 2]. One particular physical situation in which it can be gainfully employed is encountered frequently in the flow of a fluid through a porous solid. Ideally, one would like to be able to predict all quantities of interest in such a situation knowing the nature of the mixture constituents (solid and fluid in this case). But this is easier said than done, for the deformation of the solid due to the flow of fluid through it could be highly non-linear, and the fluid itself may show non-linear behaviour. Being able to accurately describe such situations quantitatively is the ultimate goal of theory of mixtures and it has been shown in [3] that starting from the balance of linear momentum, one can derive a whole gamut of mathematical models pertinent to different situations with varying levels of generality. However, using these models to actually solve realistic engineering problems is another matter altogether. Developing computational techniques to address such models is one of the important thrusts of current research. Herein, we take into account the transient and viscous effects of the fluid (while still neglecting the motion of the solid), which gives rise to *unsteady* Brinkman equation. The main intent of the present article is to present a stabilized mixed formulation for the solution of the unsteady Brinkman equation.

Darcy's equation [4] describes the flow of a fluid through a porous solid due to pressure gradients when a wide variety of assumptions are justified. It has been used to study various phenomena arising in groundwater hydrology [5, 6], enhanced oil recovery [7, 8], carbon-dioxide sequestration [9, 10] to name a few. However, it should be noted that Darcy's equation is merely an approximation to the balance of linear momentum for a fluid as it flows through a rigid porous solid (see References [3, 11, 12, 13] for a detailed derivation using the theory of interacting continua, which is also known as theory of mixtures).

Rajagopal [11] has shown that, within the context of theory of mixtures, the Darcy equation can be obtained after making a plethora of assumptions about the solid and the fluid. The assumptions include neglecting the deformations of the solid, assuming the flow to be steady, assuming a special form for the drag due to the friction at the pores, and many more. The restrictions under which the Darcy equation was derived are rather stringent, and in [11], the assumptions are systematically relaxed and an hierarchy of models is thus obtained, where the Darcy model is the one derived under the most restrictive assumptions. Moreover, it has been demonstrated in References [14, 15] that for a class of flows involving high pressure gradients, the Darcy equation is a poor approximation. However, despite its drawbacks, it remains the most popular model to describe the flow of fluids through porous solids, and for a large class of flows, it remains applicable. In the remainder of the paper, unless otherwise mentioned, the Darcy equation tacitly refers to the steady state equation.

If the viscous effects in the fluid are deemed important, as they are in many cases, and taken into account and modeled as those of the familiar Navier-Stokes fluid, the eponymous Brinkman model is obtained [16]. By further relaxing the assumption in the Brinkman model that the flow be steady, we take transient effects into account while still neglecting the inertial non-linearities, and one thus obtains (as in [11]) the *unsteady* Brinkman equation, the solution of which is the subject of the present article.

If one desires to use the finite element method to solve the Brinkman equation (or Darcy equation), then the celebrated Ladyzhenskaya-Babuška-Brezzi (LBB) condition [17, 18] must be either satisfied or circumvented. The LBB condition imposes severe restrictions on the classical mixed formulation with respect to the order of interpolation for the pressure and velocity in problems involving incompressibility. In particular, the classical mixed formulation is unstable for the equal-order interpolation for the velocity and pressure, which is computationally the most convenient. There is thus a necessity for stabilized mixed formulations, especially those that are stable under the equal-order interpolation for the velocity and pressure. Variational Multiscale (VMS) formalism provides a systematic way of developing stabilized formulation [19]. Some of the earlier works on mixed methods applied to flows through porous media are [20, 21, 22, 23, 24, 25, 26]. A thorough discussion of stabilized methods is beyond the scope of this paper, and a reader interested in this subject should refer to [27, 28, 29, 30, 30, 31, 32, 33].

**Remark 1.1.** *While it is possible to avoid a mixed formulation by using primal or single field formulations in case of the Darcy equation, such a technique cannot be applied to the* unsteady *Brinkman equation. Moreover, the primal formulation has the disadvantage of poor approximation of velocity for low-order finite elements (see* [34, 35]*).*

This VMS method for both Brinkman (or Darcy's equation) involves decomposing the velocity field into coarse/resolved scales and fine/unresolved scales. Modeling of the unresolved scales leads to a multiscale/stabilized form of the corresponding equation (see Reference [36] for details).

The subsequent attempts, quite naturally, were to use the VMS method to solve the corresponding unsteady problems. The traditional way of treating the time derivative in the unsteady equation (*unsteady* Darcy/Stokes/Navier-Stokes) has been to assume that the fine/unresolved scales are independent of time while allowing the coarse/resolved scales to depend on both the spatial and temporal variables [37, 38]. However, such an assumption is philosophically undesirable as it seems motivated by the need to somehow solve the fine scale problem rather than a sound physical basis. For completeness, we have outlined this traditional way (using a semi-discrete method) of deriving stabilized formulation in Appendix 7.2.

Another technique that has been used to treat such unsteady problems is the so called "space-time finite elements," where the temporal domain is also discretized via shape functions just like the spatial domain [39, 40, 41]. However, the disadvantage with this method is that one needs a four-dimensional mesh for a (spatially) three-dimensional problem, the implementation is not as convenient, and there are also the same issues with post-processing when visualizing the numerical results of a (spatially) three-dimensional problem.

The method followed in this article sidesteps both the preceding issues discussed in [37, 38] and [39, 40, 41]. We use Rothe's method [42] (also called the method of horizontal lines) and avoid the mathematical legerdemain involved in justifying the assumption of the time independence of the fine/unresolved scales.

The terminology "method of lines" when used with reference to numerical techniques is usually understood to imply method of "vertical" lines as distinct from the method of "horizontal" lines used in this paper. The method of lines [43, 44] is a favoured numerical technique for the solution of parabolic partial differential equations (PDEs) wherein the PDE is discretized in all but one dimension and thus reduced to a system of ordinary differential equations (ODEs). On the other hand, in the method of horizontal lines, the temporal derivative is discretized first (usually by a finite difference scheme) and the original parabolic PDE is thus reduced to a (time-discretized) elliptic PDE. The entire machinery of techniques and tools used to study elliptic PDEs (which are quite well understood) thus becomes available and therein lies the advantage of the method [45]. In a different context, the method of horizontal lines has been used before. In [46] for instance, a finite element formulation for transient incompressible viscous flows stabilized by local time-steps is derived. However, the central idea of the present paper is different; to use the method of horizontal lines in conjunction with the VMS formalism and derive a consistent formulation which circumvents an *ad hoc* assumption that is currently made.

1.1. **Main contributions of this paper.** In this paper, we systematically derive a stabilized mixed formulation for *unsteady* Brinkman equation using the method of horizontal lines and the variational multiscale formalism. The stabilization terms and parameter are obtained in a mathematically consistent manner. An important feature is that the derivation does not need the assumption that the fine-scale variables do not depend on the time, which is the case with the conventional derivation of multiscale stabilized formulations for transient mixed problems. Under the proposed formulation, equal-order interpolation for the velocity and pressure (which is computationally the most convenient) is stable. We also derive a stabilized mixed formulation for the *unsteady* Darcy equation as a special case. We also show through numerical experiments that the proposed stabilized formulation posses good spatial and temporal convergence properties.

1.2. **An outline of the paper.** In Section 2, we outline the governing equations. In Section 3, we present a stabilized mixed formulation based on the variational multiscale formalism and the method of horizontal lines for the unsteady Brinkman equation. In Section 4, the stabilized formulation obtained earlier is specialized for the case of the unsteady Darcy equation. Representative numerical results are presented in Section 5, and conclusions are drawn in Section 6.

## 2. GOVERNING EQUATIONS: UNSTEADY BRINKMAN MODEL

Let $\Omega \subset \mathbb{R}^{nd}$ be an open bounded domain, where "$nd$" is the number of spatial dimensions. Let $\Gamma$ denote the boundary, which is assumed to be smooth. A position vector is denoted by $\boldsymbol{x} \in \bar{\Omega}$, where an over-bar denotes the set closure. The time is denoted by $t \in [0, T]$, where $T$ is the total time of interest. The gradient and divergence operators with respect to $\boldsymbol{x}$ are, respectively, denoted by $\mathrm{grad}[\cdot]$ and $\mathrm{div}[\cdot]$. The velocity (vector) field is denoted as $\boldsymbol{v}(\boldsymbol{x}, t)$, and the pressure (scalar) field is denoted as $p(\boldsymbol{x}, t)$. The coefficient of viscosity, drag, and density are respectively denoted by $\mu$, $\alpha$, and $\rho$. The prescribed specific body force is denoted by $\boldsymbol{b}(\boldsymbol{x}, t)$, while the prescribed velocity on the boundary is denoted by $\boldsymbol{v}^{\mathrm{p}}(\boldsymbol{x}, t)$. The unsteady governing equations can be written as follows:

$$\rho \frac{\partial \boldsymbol{v}}{\partial t} = -\mathrm{grad}[p] - \alpha \boldsymbol{v} + \mathrm{div}[\mu \, \mathrm{grad}[\boldsymbol{v}]] + \rho \boldsymbol{b} \quad \text{in } \Omega \times (0, T) \tag{1a}$$

$$\mathrm{div}[\boldsymbol{v}] = 0 \quad \text{in } \Omega \times (0, T) \tag{1b}$$

$$\boldsymbol{v}(\boldsymbol{x}, t) = \boldsymbol{v}^{\mathrm{p}}(\boldsymbol{x}, t) \quad \text{on } \Gamma \times (0, T) \tag{1c}$$

$$\boldsymbol{v}(\boldsymbol{x}, 0) = \boldsymbol{v}_0(\boldsymbol{x}) \quad \text{in } \Omega \tag{1d}$$

where $\boldsymbol{v}_0(\boldsymbol{x})$ is the prescribed initial velocity. Few remarks are in order about the mathematical model given by equations (1a)-(1d).

**Remark 2.1.** *The drag coefficient is related to the coefficient of viscosity through* $\alpha = \dfrac{\mu}{k}$*, where* $k$ *is the coefficient of permeability. One obtains unsteady Stokes' equation if the drag term* $\alpha \boldsymbol{v}$ *is neglected. On the other hand, one obtains unsteady Darcy equation if the term* $\mathrm{div}[\mu \, \mathrm{grad}[\boldsymbol{v}]]$ *is neglected and only the normal component of the velocity on the boundary is prescribed.*

*It is noteworthy that a more general form of transient Darcy-Brinkman equation can be written as follows:*

$$\rho \frac{\partial \boldsymbol{v}}{\partial t} = -\mathrm{grad}[p] - \alpha \boldsymbol{v} + \mathrm{div}[\mu' \, \mathrm{grad}[\boldsymbol{v}]] + \rho \boldsymbol{b} \tag{2}$$

*where the factor* $\mu'$ *can be different from the coefficient of viscosity. Physically, the friction at the pores of the solid as the fluid flows past it is modeled as* $\alpha(\boldsymbol{v} - \boldsymbol{v}_{\mathrm{solid}})$*, which reduces to* $\alpha \boldsymbol{v}$ *under the assumption that the solid is rigid and does not undergo any motion. On the other hand, the viscous effects within the fluid itself are modeled as* $\mathrm{div}[\mu' \, \mathrm{grad}[\boldsymbol{v}]]$*. There are adequate discussions in the literature (including the paper by Brinkman* [47]*) about how the factor* $\mu'$ *is related to the coefficient of viscosity* $\mu$*. However, for a medium with high porosity, it is reasonable to assume that* $\mu' = \mu$ [47].

**Remark 2.2.** *It is important to note that another constraint must be enforced in order to have a unique solution for the pressure. This is usually enforced by demanding that*

$$\int_{\Omega} p \; \mathrm{d}\Omega = 0 \tag{3}$$

*In numerical simulations, however, it is more convenient to prescribe the pressure at a point, which is employed in this paper.*

**Remark 2.3.** *The balance of linear momentum for the fluid, on neglecting the non-linear convective term, takes the form*

$$\rho \frac{\partial \boldsymbol{v}}{\partial t} = \mathrm{div}[\boldsymbol{T}] + \rho \boldsymbol{b} \tag{4}$$

*Appealing to frame indifference (or even Galilean invariance), the partial stress of the fluid* $\boldsymbol{T}$ *is actually modeled as*

$$\boldsymbol{T} = -p\boldsymbol{I} + \mu \left( \mathrm{grad}[\boldsymbol{v}] + \mathrm{grad}[\boldsymbol{v}]^T \right) \tag{5}$$

*It is only when* $\mu$ *is a constant that the incompressibility constraint* $\mathrm{div}[\boldsymbol{v}] = 0$ *allows the balance of linear momentum to be simplified further and it takes the form as given in equation* (1a)*. This point becomes crucial in developing (realistic) nonlinear models in which the viscosity depends on the pressure.*

The main aim of this paper is to develop a stabilized mixed formulation for the above outlined system of partial differential equations (1a)–(1d). It is well–known that the semi-discrete formulation based on the classical mixed formulation has numerical deficiencies. In particular, the equal–order interpolation for the velocity and pressure (which is computationally the most convenient) is not stable. In the next section we develop a stabilized mixed formulation for *unsteady* Brinkman equation under which the equal–order interpolation for velocity and pressure is stable. *An important aspect is that the derivation does not need the assumption that fine-scale variables do not depend on the time, which is the case with the conventional derivation of multiscale stabilized formulations.*

## 3. A STABILIZED MIXED FORMULATION: UNSTEADY BRINKMAN EQUATION

The proposed formulation is based on the method of horizontal lines (also known as the Rothe method) [42] and variational multiscale formalism [19]. To this end, we shall discretize the time interval of interest into $N$ instants denoted by $t_n$, $(n = 0, \cdots, N)$. For simplicity, we shall assume uniform time steps, and shall denote the time step by $\Delta t := t_n - t_{n-1}$. However, it should be noted that a straightforward modification can handle non-uniform time steps. We shall denote the time discretized version of a given quantity $z(\boldsymbol{x}, t)$ at the instant of time $t_n$ as follows:

$$z^{(n)}(\boldsymbol{x}) \approx z(\boldsymbol{x}, t_n) \quad n = 0, \cdots, N \tag{6}$$

In a semi–discrete formulation, the partial differential equation (which depends on both space and time) is first spatially discretized by (say) the finite element method. This results in a system of ordinary differential equations (ODEs) or a system of differential algebraic equations (DAEs), which are numerically integrated using appropriate time stepping schemes. On the other hand, in the method of horizontal lines (also known as the Rothe's method), the partial differential equation (which depends on both space and time) is temporally discretized using a time stepping scheme. This results in another system of partial differential equations, which depends only on spatial coordinates, and the finite element method or the finite difference method is typically employed to solve this resulting system of equations. Figure 1 gives a pictorial description of the method of horizontal lines.

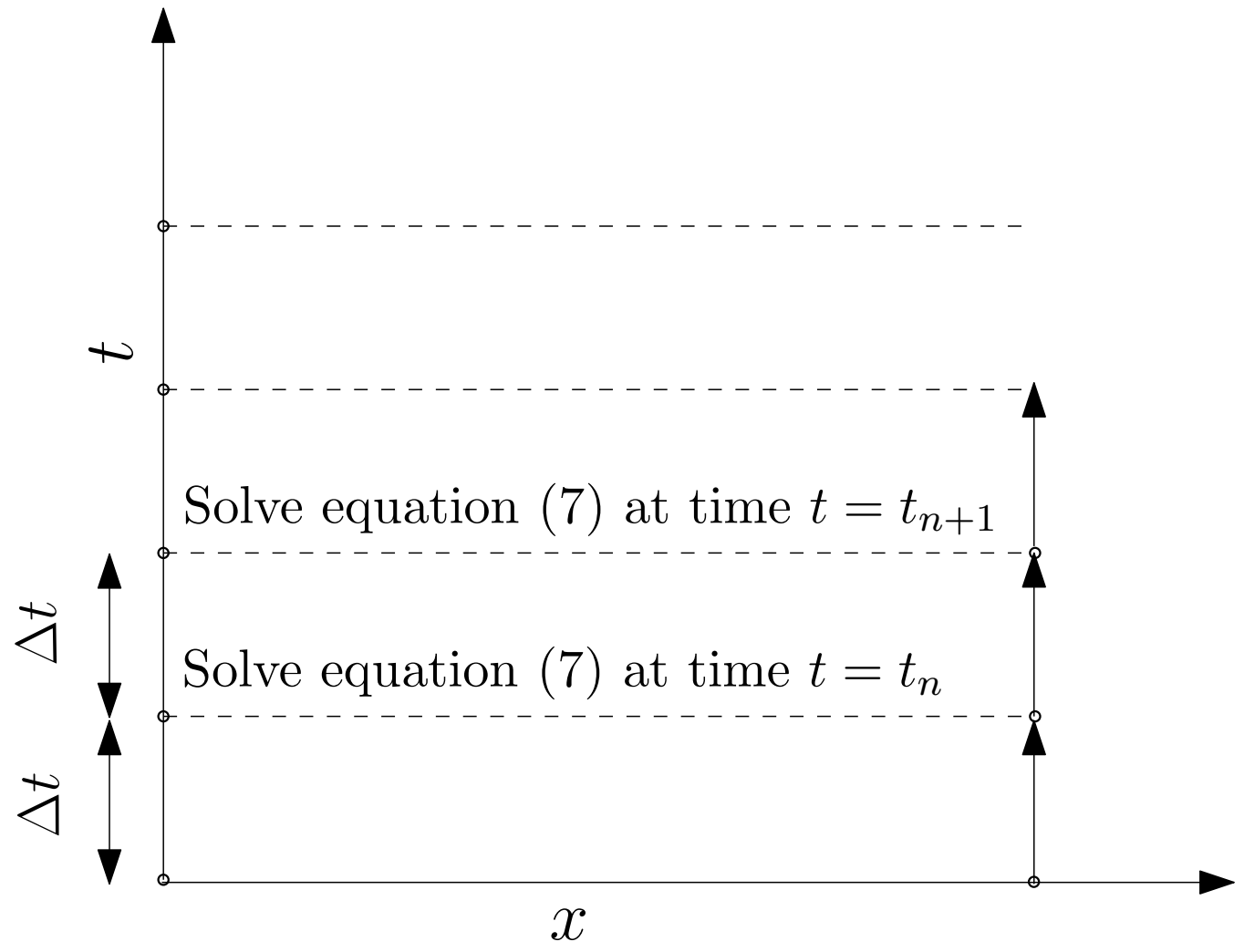

FIGURE 1. Rothe's method: The governing equations are solved at each time level after discretizing the time interval into discrete time levels.

Herein, we employ the backward Euler time stepping scheme, which is first-order accurate and unconditionally stable when applied to linear system of ordinary differential equations [48]. By employing the method of horizontal lines based on the backward Euler scheme to equation (1a) – (1d), the corresponding time discretized equations at time level $t = t_{n+1}$ can be written as follows:

$$\rho \frac{\boldsymbol{v}^{(n+1)} - \boldsymbol{v}^{(n)}}{\Delta t} = -\mathrm{grad}[p^{(n+1)}] - \alpha \boldsymbol{v}^{(n+1)} + \mathrm{div}[\mu \mathrm{grad}[\boldsymbol{v}^{(n+1)}]] + \rho \boldsymbol{b}^{(n+1)} \quad \mathrm{in}\ \Omega \tag{7a}$$

$$\mathrm{div}[\boldsymbol{v}^{(n+1)}] = 0 \quad \mathrm{in}\ \Omega \tag{7b}$$

$$\boldsymbol{v}^{(n+1)}(\boldsymbol{x}) = \boldsymbol{v}^{\mathrm{p}}(\boldsymbol{x}, t = t_{n+1}) \quad \mathrm{on}\ \Gamma \tag{7c}$$

$$\boldsymbol{v}^{(0)}(\boldsymbol{x}) = \boldsymbol{v}_0(\boldsymbol{x}) \quad \mathrm{in}\ \Omega \tag{7d}$$

**Remark 3.1.** *Equation* (7a) *can be rearranged to obtain the following equation:*

$$\left(\frac{\rho}{\Delta t} + \alpha\right) \boldsymbol{v}^{(n+1)} = -\mathrm{grad}[p^{(n+1)}] + \mathrm{div}[\mu \mathrm{grad}[\boldsymbol{v}^{(n+1)}]] + \rho \left(\boldsymbol{b}^{(n+1)} + \frac{1}{\Delta t}\boldsymbol{v}^{(n)}\right) \quad \mathrm{in}\ \Omega \tag{8}$$

*At time level $t = t_{n+1}$, the above equation* (8) *along with equations* (7b)–(7d) *resemble the well-known* steady *Brinkman equation with (modified) drag coefficient $\left(\alpha + \frac{\rho}{\Delta t}\right)$ and (modified) body force $\left(\boldsymbol{b}^{(n+1)} + \frac{1}{\Delta t}\boldsymbol{v}^{(n)}\right)$.*

For convenience, let us define $\hat{\alpha}$, $\mathcal{L}(\cdot,\cdot)$, and $\tilde{\boldsymbol{b}}^{(n+1)}(\boldsymbol{x})$ as follows:

$$\hat{\alpha} := \alpha \Delta t + \rho \tag{9a}$$

$$\mathcal{L}(\boldsymbol{w}, q) := \hat{\alpha}\boldsymbol{w} + \Delta t\ \mathrm{grad}[q] - \Delta t \mathrm{div}[\mu \mathrm{grad}[\boldsymbol{w}]] \tag{9b}$$

$$\tilde{\boldsymbol{b}}^{(n+1)} := \Delta t \boldsymbol{b}^{(n+1)} + \boldsymbol{v}^{(n)} \tag{9c}$$

Let $L^2(\Omega)$ denote the space of square integrable (scalar) functions defined on $\Omega$. The standard $L^2$ inner product over domain $K$ is denoted as $(\cdot,\cdot)_K$. That is,

$$(a; b)_K \equiv \int_K a \cdot b\ \mathrm{d}K \quad \forall a, b \in L^2(\Omega) \tag{10}$$

For simplicity, the subscript $K$ will be dropped if the domain is whole of $\Omega$ (that is, $K = \Omega$). In the remainder of the paper, we shall use the the following function spaces:

$$\mathcal{V}_t := \left\{\boldsymbol{v}(\boldsymbol{x}) \mid \boldsymbol{v}(\boldsymbol{x}) \in (H^1(\Omega))^{nd}, \boldsymbol{v}(\boldsymbol{x}) = \boldsymbol{v}^{\mathrm{p}}(\boldsymbol{x}, t)\ \mathrm{on}\ \Gamma\right\} \tag{11a}$$

$$\mathcal{W} := \left\{\boldsymbol{v}(\boldsymbol{x}) \mid \boldsymbol{v}(\boldsymbol{x}) \in (H^1(\Omega))^{nd}, \boldsymbol{v}(\boldsymbol{x}) = \boldsymbol{0}\ \mathrm{on}\ \Gamma\right\} \tag{11b}$$

$$\mathcal{P} := \left\{p(\boldsymbol{x}) \in L^2(\Omega) \mid \int_\Omega p\ \mathrm{d}\Omega = 0\right\} \tag{11c}$$

$$\mathcal{Q} := \left\{p(\boldsymbol{x}) \in H^1(\Omega) \mid \int_\Omega p\ \mathrm{d}\Omega = 0\right\} \tag{11d}$$

where $H^1(\Omega)$ is a standard Sobolev space on $\Omega$ [17]. In the next subsection, we present a derivation of the proposed stabilization mixed formulation.

**Remark 3.2.** *If traction were prescribed on a part of the boundary, then pressure would be uniquely determined and we need not enforce* $\int_\Omega p \; \mathrm{d}\Omega = 0$. *However, in this paper, we do not consider the possibility of prescribing traction.*

3.1. **A derivation of the proposed stabilized mixed formulation.** We start with the classical mixed formulation for equations (7a)–(7d), which can written as: Find $\boldsymbol{v}^{(n+1)}(\boldsymbol{x}) \in \mathcal{V}_{t=t_{n+1}}$ and $p^{(n+1)}(\boldsymbol{x}) \in \mathcal{P}$ such that we have

$$
\begin{aligned}
&(\boldsymbol{w}; \hat{\alpha}\boldsymbol{v}^{(n+1)}) + \Delta t(\mathrm{grad}[\boldsymbol{w}]; \mu \; \mathrm{grad}[\boldsymbol{v}^{(n+1)}]) - \Delta t(\mathrm{div}[\boldsymbol{w}]; p^{(n+1)}) - \Delta t(q; \mathrm{div}[\boldsymbol{v}^{(n+1)}]) \\
&\quad = (\boldsymbol{w}; \rho \tilde{\boldsymbol{b}}^{(n+1)}) \quad \forall \boldsymbol{w}(\boldsymbol{x}) \in \mathcal{W}, \; q(\boldsymbol{x}) \in \mathcal{P} \qquad (12)
\end{aligned}
$$

As mentioned earlier, the classical mixed formulation has numerical instabilities, and hence the need for alternate (stabilized) formulations. To this end, let us divide the domain $\Omega$ into "$Nele$" non-overlapping subdomains $\Omega^e$ (which, in the finite element context, will be elements) such that

$$
\bar{\Omega} = \bigcup_{e=1}^{Nele} \bar{\Omega}^e \qquad (13)
$$

The boundary of the element $\Omega^e$ is denoted by $\Gamma^e := \bar{\Omega}^e - \Omega^e$. We decompose the velocity field $\boldsymbol{v}^{(n+1)}(\boldsymbol{x})$ into coarse-scale and fine-scale components, indicated as $\boldsymbol{v}_c^{(n+1)}(\boldsymbol{x})$ and $\boldsymbol{v}_f^{(n+1)}(\boldsymbol{x})$, respectively. That is,

$$
\boldsymbol{v}^{(n+1)}(\boldsymbol{x}) = \boldsymbol{v}_c^{(n+1)}(\boldsymbol{x}) + \boldsymbol{v}_f^{(n+1)}(\boldsymbol{x}) \qquad (14)
$$

Likewise, we decompose the weighting function $\boldsymbol{w}(\boldsymbol{x})$ into coarse-scale $\boldsymbol{w}_c(\boldsymbol{x})$ and fine-scale $\boldsymbol{w}_f(\boldsymbol{x})$ components. At the outset, we recognize that this multi-scale decomposition can be done irrespective of whether the problem is linear or non-linear. We further make an assumption that the fine-scale components vanish along each element boundary. That is,

$$
\boldsymbol{v}_f^{(n+1)}(\boldsymbol{x}) = \boldsymbol{w}_f(\boldsymbol{x}) = \boldsymbol{0} \quad \text{on} \quad \Gamma^e \; ; \; e = 1, \cdots, Nele \qquad (15)
$$

Let $\mathcal{V}_c$ be the function space for the coarse-scale component of the velocity $\boldsymbol{v}_c^{(n+1)}(\boldsymbol{x})$, and $\mathcal{W}_c$ be the function space for $\boldsymbol{w}_c(\boldsymbol{x})$. The spaces $\mathcal{V}_c$ and $\mathcal{W}_c$ are defined as

$$
\mathcal{V}_c := \mathcal{V}_t; \; \mathcal{W}_c := \mathcal{W} \qquad (16)
$$

where $\mathcal{V}_t$ and $\mathcal{W}$ are as defined earlier in equation (11a) and equation (11b), respectively. Let $\mathcal{V}_f$ be the function space for both the fine-scale component of the velocity $\boldsymbol{v}_f^{(n+1)}(\boldsymbol{x})$ and its corresponding weighting function $\boldsymbol{w}_f(\boldsymbol{x})$. It is defined as follows:

$$
\mathcal{V}_f := \{\boldsymbol{v}(\boldsymbol{x}) \in \left(H^1(\Omega)\right)^{nd} \mid \boldsymbol{v}(\boldsymbol{x}) = \boldsymbol{0} \text{ on } \Gamma^e; e = 1, \cdots, Nele\} \qquad (17)
$$

Substituting equation (14) into equation (12), one obtains two subproblems. The coarse-scale and fine-scale subproblems are, respectively, defined as follows: Find $\boldsymbol{v}_c^{(n+1)}(\boldsymbol{x}) \in \mathcal{V}_c$, $\boldsymbol{v}_f^{(n+1)}(\boldsymbol{x}) \in \mathcal{V}_f$, and $p^{(n+1)}(\boldsymbol{x}) \in \mathcal{P}$ such that we have

$$\begin{aligned}
&\left(\boldsymbol{w}_c; \hat{\alpha}\boldsymbol{v}_c^{(n+1)}\right) + (\boldsymbol{w}_c; \hat{\alpha}\boldsymbol{v}_f^{(n+1)}) + \Delta t \left(\mathrm{grad}[\boldsymbol{w}_c]; \mu \; \mathrm{grad}[\boldsymbol{v}_c^{(n+1)}]\right) \\
&+ \Delta t(\mathrm{grad}[\boldsymbol{w}_c]; \mu \; \mathrm{grad}[\boldsymbol{v}_f^{(n+1)}]) - \Delta t \left(\mathrm{div}[\boldsymbol{w}_c]; p^{(n+1)}\right) - \Delta t \left(q; \mathrm{div}[\boldsymbol{v}_c^{(n+1)}]\right) \\
&- \Delta t(q; \mathrm{div}[\boldsymbol{v}_f^{(n+1)}]) = (\boldsymbol{w}_c; \rho\tilde{\boldsymbol{b}}^{(n+1)}) \quad \forall \boldsymbol{w}_c(\boldsymbol{x}) \in \mathcal{W}_c, \; q(\boldsymbol{x}) \in \mathcal{P} \qquad \text{(18a)} \\
&(\boldsymbol{w}_f; \hat{\alpha}\boldsymbol{v}_c^{(n+1)}) + (\boldsymbol{w}_f; \hat{\alpha}\boldsymbol{v}_f^{(n+1)}) + \Delta t \left(\mathrm{grad}[\boldsymbol{w}_f]; \mu \; \mathrm{grad}[\boldsymbol{v}_c^{(n+1)}]\right) \\
&+ \Delta t(\mathrm{grad}[\boldsymbol{w}_f]; \mu \; \mathrm{grad}[\boldsymbol{v}_f^{(n+1)}]) - \Delta t(\mathrm{div}[\boldsymbol{w}_f]; p^{(n+1)}) = (\boldsymbol{w}_f; \rho\tilde{\boldsymbol{b}}^{(n+1)}) \quad \forall \boldsymbol{w}_f(\boldsymbol{x}) \in \mathcal{W}_f \qquad \text{(18b)}
\end{aligned}$$

The fine scale problem (18b) can be solved (either approximately or exactly) and substituted into the coarse scale problem (18a) to eliminate the fine scale variables completely. But before we attempt to do that, it is necessary to interpolate the fine-scale variables. [We interpolate the coarse-scale variables using the standard Lagrangian shape functions.] The fine-scale variables are interpolated using bubble functions. A bubble function defined over an element is a function that vanishes on the boundary of the element. The fine-scale variables over an element are interpolated as follows:

$$\boldsymbol{v}_f^{(n+1)}(\boldsymbol{x}) = b^e(\boldsymbol{x})\boldsymbol{\beta}, \quad \boldsymbol{w}_f(\boldsymbol{x}) = b^e(\boldsymbol{x})\boldsymbol{\gamma} \qquad (19)$$

where $b^e(\boldsymbol{x})$ is the bubble function while $\boldsymbol{\beta}$ and $\boldsymbol{\gamma}$ are column vectors of size $nd \times 1$. It will be easier to construct a bubble function on the reference element, and use the isoparametric mapping to define it with respect to the original coordinate system. Therefore, the bubble functions will be functions of $\boldsymbol{\zeta}$, which is the position vector in isoparametric coordinates. An illustration of all the different bubbles used in this paper is provided in Figure 2.

Solving (18b) for $\boldsymbol{\beta}$ and substituting into (18a), the proposed stabilized mixed formulation at time level $t = t_{n+1}$ we arrive at is the following (dropping the suffix $c$ for convenience): Find $\boldsymbol{v}^{(n+1)}(\boldsymbol{x}) \in \mathcal{V}_{t=t_{n+1}}$ and $p^{(n+1)}(\boldsymbol{x}) \in \mathcal{Q}$ such that we have

$$\begin{aligned}
&(\boldsymbol{w}; \hat{\alpha}\boldsymbol{v}^{(n+1)}) + \Delta t(\mathrm{grad}[\boldsymbol{w}]; \mu \; \mathrm{grad}[\boldsymbol{v}^{(n+1)}]) - \Delta t(\mathrm{div}[\boldsymbol{w}]; p^{(n+1)}) - \Delta t(q; \mathrm{div}[\boldsymbol{v}^{(n+1)}]) \\
&- \sum_{e=1}^{Nele} \left(\mathcal{L}(\boldsymbol{w}, q); \tau(\boldsymbol{x})\mathcal{L}(\boldsymbol{v}^{(n+1)}, p^{(n+1)})\right)_{\Omega^e} = (\boldsymbol{w}; \rho\tilde{\boldsymbol{b}}^{(n+1)}) - \sum_{e=1}^{Nele} \left(\mathcal{L}(\boldsymbol{w}, q); \tau(\boldsymbol{x})\rho\tilde{\boldsymbol{b}}^{(n+1)}\right)_{\Omega^e} \\
&\forall \boldsymbol{w}(\boldsymbol{x}) \in \mathcal{W}, \; q(\boldsymbol{x}) \in \mathcal{Q} \qquad (20)
\end{aligned}$$

where the stabilization parameter $\tau(\boldsymbol{x})$ is given by

$$\tau(\boldsymbol{x}) = b^e(\boldsymbol{x}) \left(\int_{\Omega^e} b^e(\boldsymbol{x}) \; \mathrm{d}\Omega\right) \left(\int_{\Omega^e} \left(\mu \; \Delta t \; \|\mathrm{grad}[b^e]\|^2 + \hat{\alpha}(b^e)^2\right) \; \mathrm{d}\Omega\right)^{-1} \qquad (21)$$

where $b^e(\boldsymbol{x})$ is a bubble function. Thus we obtain a formulation which has additional stabilizing terms compared to classical mixed formulation. A systematic numerical implementation of the proposed formulation is outlined in Algorithm 1. A few remarks about the stabilization parameter are in order.

**Remark 3.3.** *We have $\hat{\alpha} > 0,\ \mu > 0,\ \Delta t > 0$. Also, for $\boldsymbol{x} \in \Omega^e,\ (b^e(\boldsymbol{x}))^2 > 0, \|\mathrm{grad}[b^e]\|^2 > 0$. The stabilization parameter $\tau(\boldsymbol{x})$ is thus well-defined for all permissible $\hat{\alpha},\ \mu$ and $\Delta t$ since the integrand $\left(\mu\, \Delta t\, \|\mathrm{grad}[b^e]\|^2 + \hat{\alpha}(b^e)^2\right) > 0$.*

**Remark 3.4.** *On defining*

$$\tau_{\mathrm{avg}} := \frac{1}{\mathrm{meas}(\Omega^e)} \int_{\Omega^e} \tau(\boldsymbol{x})\ \mathrm{d}\Omega \tag{22}$$

*one can get bounds on $\tau_{\mathrm{avg}}$ for some special cases of interest. The Poincaré inequality* [49] *implies that*

$$\int_{\Omega^e} (b^e(\boldsymbol{x}))^2\ \mathrm{d}\Omega \leq \gamma(\Omega^e) \int_{\Omega^e} \|\mathrm{grad}[b^e]\|^2\ \mathrm{d}\Omega \tag{23}$$

*where $\gamma(\Omega^e)$ is a domain-dependent positive constant. Note that, in the above equation, we have used the fact that the bubble function vanishes on the boundary. In particular, we have assumed that $b^e(\boldsymbol{x}) \in H^1_0(\Omega^e)$. The elements in a typical finite element mesh will all be convex domains. For such domains, there are good estimates for the parameter $\gamma$ (for example, see Reference* [50]*). The Cauchy-Schwarz inequality on the other hand implies*

$$\left(\int_{\Omega^e} b^e(\boldsymbol{x})\ \mathrm{d}\Omega\right)^2 \leq \mathrm{meas}(\Omega^e) \left(\int_{\Omega^e} (b^e(\boldsymbol{x}))^2\ \mathrm{d}\Omega\right) \tag{24}$$

*The inequalities* (23) *and* (24) *immediately yield*

$$0 < \tau_{\mathrm{avg}} \leq (\hat{\alpha} + \mu\, \Delta t\, \gamma(\Omega^e))^{-1} \tag{25}$$

*Specifically, for the unsteady Darcy equation ($\mu = 0$) and the unsteady Stokes equation ($\hat{\alpha} = \rho$), we have*

$$0 < \tau_{\mathrm{avg}} \leq (\hat{\alpha})^{-1} \quad \text{if } \mu = 0 \tag{26a}$$

$$0 < \tau_{\mathrm{avg}} \leq (\rho + \mu\, \Delta t\, \gamma(\Omega^e))^{-1} \quad \text{if } \hat{\alpha} = \rho \tag{26b}$$

**Remark 3.5.** *Readers consulting the references provided herein are forewarned of a few errors and oversights in* [51] *. The boundary condition prescribed in equation (3) of* [51] *is appropriate to the Darcy equation and not to the Brinkman equation. A similar comment applies to equation (4) where the function space for the velocity is defined; such a space is valid for Darcy equation but not for the Brinkman equation. Moreover, in equation (4), it ought to read as* $\mathrm{trace}(\boldsymbol{v}\cdot\boldsymbol{n}) = \psi \in H^{-1/2}(\Gamma)$.

*Finally, the stabilization parameter defined in equation (30) will make sense only if we understand* $\mu\Delta\boldsymbol{v}$ *to mean* $\mu\Delta\boldsymbol{v} + \mu\mathrm{div}[\mathrm{grad}[\boldsymbol{v}]^T]$.

**Algorithm 1** Implementation of the proposed formulation

Input: Initial condition $\boldsymbol{v}_0$, Time period of integration $T$, time step $\Delta t$
Set $\boldsymbol{v}^{(n)} = \boldsymbol{v}_0$
Set $t = 0$
**while** $t < T$ **do**
    $\Delta t = \min(\Delta t, T - t), \quad t = t + \Delta t$
    Using $\boldsymbol{v}^{(n)}$, solve equation (20) to obtain $\boldsymbol{v}^{(n+1)}, p^{(n+1)}$
    Set $\boldsymbol{v}^{(n)} = \boldsymbol{v}^{(n+1)}$
**end while**

**Remark 3.6.** *In the implementation of the algorithm, all the terms in equation* (20) *need not be evaluated repeatedly at each time step since most of them do not depend on the temporal variable. In fact, only the terms involving* $\tilde{\boldsymbol{b}}$ *need to be evaluated repeatedly.*

## 4. SPECIAL CASE : A STABILIZED MIXED FORMULATION FOR UNSTEADY DARCY EQUATION

We now present a stabilized mixed formulation for the unsteady Darcy equation following the same procedure as outlined in the previous section. Most of the quantities appearing have been defined previously. The ones appearing for the first time are defined here. The boundary $\Gamma$ is divided into two parts, denoted by $\Gamma^v$ and $\Gamma^p$, such that $\Gamma^v \cap \Gamma^p = \emptyset$ and $\Gamma^v \cup \Gamma^p = \Gamma$. $\Gamma^p$ is the part of the boundary on which the pressure $p$ is prescribed, and $\Gamma^v$ is that part of the boundary on which the normal component of the velocity is prescribed. The prescribed normal component of the velocity on the boundary is denoted by $\psi(\boldsymbol{x}, t)$, and the prescribed pressure is denoted by $p_0(\boldsymbol{x}, t)$. The unsteady governing equations can be written as

$$\rho \frac{\partial \boldsymbol{v}}{\partial t} = -\mathrm{grad}[p] - \alpha \boldsymbol{v} + \rho \boldsymbol{b} \quad \text{in } \Omega \times (0, T) \tag{27a}$$

$$\mathrm{div}[\boldsymbol{v}] = 0 \quad \text{in } \Omega \times (0, T) \tag{27b}$$

$$p(\boldsymbol{x}, t) = p_0(\boldsymbol{x}, t) \quad \text{on } \Gamma^p \times (0, T) \tag{27c}$$

$$\boldsymbol{v}(\boldsymbol{x}, t) \cdot \boldsymbol{n}(\boldsymbol{x}) = \psi(\boldsymbol{x}, t) \quad \text{on } \Gamma^v \times (0, T) \tag{27d}$$

$$\boldsymbol{v}(\boldsymbol{x}, 0) = \boldsymbol{v}_0(\boldsymbol{x}) \quad \text{in } \Omega \tag{27e}$$

where $\boldsymbol{n}(\boldsymbol{x})$ is the unit outward normal on the boundary.

**Remark 4.1.** *It is important to note that if $\Gamma^v = \Gamma$, then for a given $\Omega$ (and hence given $\Gamma$) $\psi$ cannot be arbitrarily specified. Using the divergence theorem it is apparent that $\psi$ must meet the following compatibility condition:*

$$\int_{\Gamma=\Gamma^v} \psi \; \mathrm{d}\Gamma = 0 \tag{28}$$

The relevant function spaces for the velocity and pressure fields in case of the unsteady Darcy equation are defined as follows:

$$\bar{\mathcal{V}}_t := \left\{ \boldsymbol{v}(\boldsymbol{x}) \mid \boldsymbol{v}(\boldsymbol{x}) \in H^1(\mathrm{div};\Omega), \boldsymbol{v}(\boldsymbol{x}) \cdot \boldsymbol{n}(\boldsymbol{x}) = \psi(\boldsymbol{x},t) \text{ on } \Gamma^v \right\} \tag{29a}$$

$$\bar{\mathcal{W}} := \left\{ \boldsymbol{v}(\boldsymbol{x}) \mid \boldsymbol{v}(\boldsymbol{x}) \in H^1(\mathrm{div};\Omega), \boldsymbol{v}(\boldsymbol{x}) \cdot \boldsymbol{n}(\boldsymbol{x}) = 0 \text{ on } \Gamma^v \right\} \tag{29b}$$

$$\bar{\mathcal{P}} := \left\{ p(\boldsymbol{x}) \mid p(\boldsymbol{x}) \in L^2(\Omega) \right\}, \quad \bar{\mathcal{Q}} := \left\{ p(\boldsymbol{x}) \mid p(\boldsymbol{x}) \in H^1(\Omega) \right\} \tag{29c}$$

where $H^1(\mathrm{div};\Omega)$ is defined by

$$H^1(\mathrm{div};\Omega) := \left\{ \boldsymbol{v}(\boldsymbol{x}) \in \left(L^2(\Omega)\right)^{nd} \mid \mathrm{div}[\boldsymbol{v}] \in L^2(\Omega) \right\} \tag{30}$$

**Remark 4.2.** *If we have $\Gamma^v = \Gamma$, then*

$$\bar{\mathcal{P}} := \left\{ p(\boldsymbol{x}) \in L^2(\Omega) \mid \int_\Omega p \; \mathrm{d}\Omega = 0 \right\}, \quad \bar{\mathcal{Q}} := \left\{ p(\boldsymbol{x}) \in H^1(\Omega) \mid \int_\Omega p \; \mathrm{d}\Omega = 0 \right\} \tag{31}$$

For convenience, let us define the operator $\bar{\mathcal{L}}(\cdot,\cdot)$, $\hat{\alpha}$ and $\tilde{\boldsymbol{b}}^{(n+1)}(\boldsymbol{x})$ to be

$$\bar{\mathcal{L}}(\boldsymbol{w}, q) := \hat{\alpha}\boldsymbol{w} + \Delta t \; \mathrm{grad}[q] \tag{32a}$$

$$\hat{\alpha} := \rho + \alpha \Delta t \tag{32b}$$

$$\tilde{\boldsymbol{b}}^{(n+1)} := \Delta t \boldsymbol{b}^{(n+1)} + \boldsymbol{v}^{(n)} \tag{32c}$$

A stabilized mixed formulation for unsteady Darcy equation at time level $t = t_{n+1}$ can be written as: Find $\boldsymbol{v}^{(n+1)}(\boldsymbol{x}) \in \bar{\mathcal{V}}_{t=t_{n+1}}$ and $p^{(n+1)}(\boldsymbol{x}) \in \bar{\mathcal{Q}}$ such that we have

$$\begin{aligned} &(\boldsymbol{w}; \hat{\alpha}\boldsymbol{v}^{(n+1)}) - \Delta t(\mathrm{div}[\boldsymbol{w}]; p^{(n+1)}) - \Delta t(q; \mathrm{div}[\boldsymbol{v}^{(n+1)}]) \\ &\quad - \left( \bar{\mathcal{L}}(\boldsymbol{w}, q); \tau(\boldsymbol{x})\hat{\alpha}^{-1} \bar{\mathcal{L}}(\boldsymbol{v}^{(n+1)}, p^{(n+1)}) \right) \\ &\quad = (\boldsymbol{w}; \rho\tilde{\boldsymbol{b}}^{(n+1)}) - \Delta t(\boldsymbol{w} \cdot \boldsymbol{n}; {p_0}^{(n+1)})_{\Gamma^p} \\ &\quad - \left( \bar{\mathcal{L}}(\boldsymbol{w}, q); \tau(\boldsymbol{x})\hat{\alpha}^{-1} \rho \tilde{\boldsymbol{b}}^{(n+1)} \right) \quad \forall \boldsymbol{w}(\boldsymbol{x}) \in \mathcal{W}, \; q(\boldsymbol{x}) \in \mathcal{Q} \end{aligned} \tag{33}$$

where the stabilization parameter $\tau(\boldsymbol{x})$ is given by

$$\tau(\boldsymbol{x}) = b^e(\boldsymbol{x}) \left[ \int_{\Omega^e} b^e(\boldsymbol{x}) \; \mathrm{d}\Omega \right] \left[ \int_{\Omega^e} (b^e)^2 \; \mathrm{d}\Omega \right]^{-1} \tag{34}$$

For simplicity, one can use a constant representative value of $\tau(\boldsymbol{x}) = \tau_{\mathrm{avg}} = \dfrac{1}{2}$ for all elements. A careful discussion on the use of the representative value for $\tau_{\mathrm{avg}}$ can be found in References [36, 35]).

**Remark 4.3.** *While this formulation requires the pressure $p \in H^1(\Omega)$, it is possible to use a classical mixed formulation with different interpolations for the velocity and pressure so that we only require $p \in L^2(\Omega)$. But the relaxation in the requirement comes at the cost of increasing complexity and greater difficulty in implementation. Moreover, the advantage is short-lived because as we proceed in the hierarchy to a Brinkman model, we will require further smoothness from the pressure field in any case.*

## 5. REPRESENTATIVE NUMERICAL RESULTS

In this section, we shall show the performance of the proposed stabilized formulation using representative numerical examples. We also perform spatial and temporal numerical convergence studies. In all the numerical simulations, we have employed the equal-order interpolation for the velocity and pressure. The drag coefficient, density, and viscosity are all taken as $\alpha = \rho = \mu = 1$ for all the cases.

# Part 1. UNSTEADY DARCY EQUATION

### 5.1. One-dimensional problem.

Consider the domain to be of unit length. The body force, boundary conditions, and initial conditions are, respectively, given by

$$b(x,t) = -\sin(t) + -4x^3 \sin(t) + \cos(t) \tag{35a}$$

$$v(0,t) = v(1,t) = \cos(t) \tag{35b}$$

$$v(x,0) = 1 \tag{35c}$$

For uniqueness, the pressure is prescribed at $x = 0$, and the value is taken to be zero. The analytical solution is given by

$$v(x,t) = \cos(t) \tag{36a}$$

$$p(x,t) = -x^4 \sin(t) \tag{36b}$$

The above problem is solved using the proposed formulation, and the computational domain is divided into 20 equal-sized two-node linear finite elements. In Figures 3 and 4 the obtained numerical solution is compared with the exact solution for time instants at $T = 1$ and $T = 2$, respectively.

#### 5.1.1. *Temporal numerical convergence studies.*

To study the temporal convergence, we compared the numerical solutions obtained at time $T = 1$ using successive time steps of magnitude $10^{-2}$, $\frac{10^{-2}}{2}$, $\frac{10^{-2}}{4}$, $\frac{10^{-2}}{8}$ and $\frac{10^{-2}}{16}$ with a mesh of 128 elements. The terminal rate of convergence is found to be approximately first order (see Figure 5), which is what is expected of a backward Euler time discretization scheme.

5.1.2. *Spatial numerical convergence studies.* To study the spatial convergence, we compared the numerical solutions obtained at time $T = 1$ using hierarchical meshes (consisting of 8, 16, 32, 64 and 128 elements) with a time step of $10^{-4}$. The terminal rates of convergence are found to be (approximately) first order for the pressure and second order for the velocity, which are shown Figure 5.

5.2. **Two-dimensional problem.** The computational domain is a bi-unit square. The initial condition for the velocity is taken to be zero, and homogeneous boundary conditions for the normal component of the velocity are prescribed on the whole boundary. The test problem is pictorially defined in Figure 6. The body force is given by

$$b_x(x, y, t) = (2\exp(t) - 1)\sin(2\pi x)\sin(2\pi y) + \sin(\pi t)y(1 - y)(1 - 2x) \tag{37a}$$

$$b_y(x, y, t) = (2\exp(t) - 1)(\cos(2\pi x)\cos(2\pi y) - \cos(2\pi x)) + \sin(\pi t)x(1 - x)(1 - 2y) \tag{37b}$$

The analytical solution is given by

$$v_x(x, y, t) = (\exp(t) - 1)\sin(2\pi x)\sin(2\pi y), \tag{38a}$$

$$v_y(x, y, t) = (\exp(t) - 1)(\cos(2\pi x)\cos(2\pi y) - \cos(2\pi x)) \tag{38b}$$

$$p(x, y, t) = x(1 - x)y(1 - y)\sin(\pi t) \tag{38c}$$

This problem was solved in turn on a T3, Q4 and an unstructured T3 grid respectively (see Figure 7). A grid of size $31 \times 31$ was used for the case of T3 and Q4 elements while 448 elements were used in the case of the unstructured T3 grid. A comparison of the exact and numerical solution at time $T = 0.5$ is presented in the contour plots in Figures 8 – 10.

5.2.1. *Spatial convergence.* The spatial convergence was studied by comparing numerical solutions obtained at time $T = 0.5$ using successively finer meshes with 5, 9, 17 and 33 nodes along each side respectively with a time step of magnitude $10^{-3}$. The terminal rate of convergence is found to be approximately second order for both pressure and velocity with both four-node quadrilateral and three-node triangular finite elements (see Figure 11).

## Part 2. UNSTEADY BRINKMAN EQUATION

5.3. **One-dimensional problem.** Consider the domain to be of unit length. The body force, boundary conditions, and initial conditions are, respectively, given by

$$b(x, t) = \cos(10\pi t) + 6\pi\exp(t)\cos(6\pi x) - 10\pi\sin(10\pi t) \tag{39a}$$

$$v(0, t) = v(1, t) = \cos(10\pi t) \tag{39b}$$

$$v(x, 0) = 1 \tag{39c}$$

For uniqueness, the pressure is prescribed at $x = 0$, and the value is taken to be zero. The analytical solution is given by

$$v(x,t) = \cos(10\pi t) \tag{40a}$$

$$p(x,t) = \exp(t)\sin(6\pi x) \tag{40b}$$

The above problem is solved using the proposed formulation, and the computational domain is divided into 40 equal-sized two-node linear finite elements. In Figures 12 and 13 the obtained numericalsolution is compared with the exact solution for time instants at $T = 0.1$ and $T = 0.2$, respectively.

5.3.1. *Temporal numerical convergence studies.* To study the temporal convergence, we compared the numerical solutions obtained at time $T = 1$ using successive time steps of magnitude $\dfrac{10^{-2}}{2}$, $\dfrac{10^{-2}}{4}$, $\dfrac{10^{-2}}{8}$ and $\dfrac{10^{-2}}{16}$ with a mesh of 128 elements. The terminal rate of convergence is found to be approximately first order (see Figure 14), which is what is expected of a backward Euler time discretization scheme.

5.3.2. *Spatial numerical convergence studies.* To study the spatial convergence, we compared the numerical solutions obtained at time $T = 0.01$ using hierarchical meshes (consisting of 16, 32, 64, 128 and 256 elements) with a time step of $10^{-5}$. The terminal rates of convergence are found to be (approximately) first order for the pressure and second order for the velocity, which are shown Figure 14.

5.4. **Two-dimensional problem.** The computational domain is a bi-unit square. The initial condition for the velocity is taken to be zero, and velocity is prescribed everywhere on the boundary. The test problem is pictorially defined in Figure . The body force is given by

$$\begin{aligned} b_x(x,y,t) &= (2\exp(t) - 1)\sin(2\pi x)\sin(2\pi y) + 8\pi^2(\exp(t) - 1)\sin(2\pi x)\sin(2\pi y) \\ &\quad + \pi\exp(t)\cos(\pi x)\sin(\pi y) \end{aligned} \tag{41a}$$

$$\begin{aligned} b_y(x,y,t) &= (2\exp(t) - 1)(\cos(2\pi x)\cos(2\pi y) - \cos(2\pi x)) + 8\pi^2(\exp(t) - 1)\cos(2\pi x)\cos(2\pi y) \\ &\quad - 4\pi^2(\exp(t) - 1)\cos(2\pi x) + \pi\exp(t)\cos(\pi y)\sin(\pi x) \end{aligned} \tag{41b}$$

The analytical solution is given by

$$v_x(x,y,t) = (\exp(t) - 1)\sin(2\pi x)\sin(2\pi y), \tag{42a}$$

$$v_y(x,y,t) = (\exp(t) - 1)(\cos(2\pi x)\cos(2\pi y) - \cos(2\pi x)) \tag{42b}$$

$$p(x,y,t) = \sin(\pi x)\sin(\pi y)\exp(t); \tag{42c}$$

This problem was solved in turn on a (six-node triangle element) T6 and (nine-node quadrilateral element) Q9 grid respectively (see Figure 2). A grid of size $31 \times 31$ was used for the case of Q9 while a $41 \times 41$ grid was used in the case of the T6 element. A comparison of the exact and numerical solution at time $T = 0.20$ is presented in the contour plots in Figures 15 and 16.

**Remark 5.1.** *The T3 and Q4 elements perform satisfactorily for the unsteady Darcy equation. However in case of the unsteady Brinkman equation, they perform poorly with respect to computation of the pressure, mainly because the second order term (Laplacian) is identically zero. Hence we have used higher order elements, namely the T6 and Q9 elements for the unsteady Brinkman equation. However, with higher order elements, the choice of bubble functions is quite critical and is by itself an important question. In this case, with better choice of bubbles, the results should improve.*

5.4.1. *Spatial convergence.* The spatial convergence was studied at two different time instants, $T = 0.01$ and $T = 0.1$, by comparing numerical solutions obtained using successively finer meshes with a constant time step. For the case of $T = 0.01$, we use a time step $\Delta t = 10^{-4}$ while for $T = 0.1$, $\Delta t = 10^{-3}$. In case of the T6 element, meshes consisting of 5, 9, 17 and 33 nodes along each side were used in succession while for the Q9 element, in turn, it consisted of 3, 5, 9, and 17 nodes along each side. The terminal rate of convergence is found to be approximately second order for both pressure and velocity with both Q9 and T6 finite elements (see Figure 17 and 18).

## 6. CONCLUSIONS

In this paper, we have presented a stabilized mixed formulation for *unsteady* Brinkman equation. The proposed formulation is consistently derived based on the variational multiscale formalism combined with the method of horizontal lines. The derivation does not require the assumption that the fine-scale terms do not depend on time, which is required under many of the existing stabilized variational multiscale formulations for transient problems. Under the proposed formulation, it has been shown numerically that the equal-order interpolation for the velocity and pressure is stable (which is not the case with the classical mixed formulation). The performance of the proposed formulation is illustrated using various representative numerical examples. Spatial and temporal numerical convergence studies are also performed, and the proposed formulation performed well.

## 7. APPENDIX

In Appendix, we present a finite element discretization of $\Delta \boldsymbol{v}$ and $\Delta \boldsymbol{w}$. For comparison and completeness, we also present the conventional way of deriving stabilized mixed formulation based on the variational multiscale approach.

7.1. **Finite element discretization of $\Delta \boldsymbol{v}$ and $\Delta \boldsymbol{w}$ terms.** We first introduce some notation and definitions, which be useful in writing the "stiffness" matrices in a systematic manner. Let $\boldsymbol{A}$ and $\boldsymbol{B}$ be matrices of size $n \times m$ and $p \times q$, respectively.

$$\boldsymbol{A} = \begin{bmatrix} a_{1,1} & \dots & a_{1,m} \\ \vdots & \ddots & \vdots \\ a_{n,1} & \dots & a_{n,m} \end{bmatrix} ; \ \boldsymbol{B} = \begin{bmatrix} b_{1,1} & \dots & b_{1,q} \\ \vdots & \ddots & \vdots \\ b_{p,1} & \dots & b_{p,q} \end{bmatrix}$$

The *Kronecker product* of these matrices is an $np \times mq$ matrix, and is defined as

$$\boldsymbol{A} \odot \boldsymbol{B} := \begin{bmatrix} a_{1,1}\boldsymbol{B} & \dots & a_{1,m}\boldsymbol{B} \\ \vdots & \ddots & \vdots \\ a_{n,1}\boldsymbol{B} & \dots & a_{n,m}\boldsymbol{B} \end{bmatrix}$$

Note that the Kronecker product can defined for *any* two given matrices (irrespective of their dimensions). The vec[$\cdot$] operator is defined as

$$\mathrm{vec}[\boldsymbol{A}] := \begin{bmatrix} a_{1,1} \\ \vdots \\ a_{1,m} \\ \vdots \\ a_{n,1} \\ \vdots \\ a_{n,m} \end{bmatrix}$$

Some relevant properties of Kronecker product and vec[$\cdot$] operator are as follows.

- $\mathrm{vec}[\boldsymbol{ACB}] = \left(\boldsymbol{B}^T \odot \boldsymbol{A}\right) \mathrm{vec}[\boldsymbol{C}]$
- $(\boldsymbol{A} \odot \boldsymbol{B})\,(\boldsymbol{C} \odot \boldsymbol{D}) = (\boldsymbol{AC} \odot \boldsymbol{BD})$
- $\mathrm{vec}[\boldsymbol{A} + \boldsymbol{B}] = \mathrm{vec}[\boldsymbol{A}] + \mathrm{vec}[\boldsymbol{B}]$

For convenience and effective computer implementation, we shall represent a three-dimensional array as a two-dimensional matrix. Consider a three-dimensional array $\boldsymbol{H}$ of size $m \times n \times p$, and its corresponding matrix representation is denoted by $\mathrm{mat}_1(\cdot)$ and is defined as follows:

$$\mathrm{mat}_1[\boldsymbol{H}] := \begin{pmatrix} H_{111} & \cdots & H_{1n1} & \cdots & \cdots & H_{11p} & \cdots & H_{1np} \\ \vdots & \ddots & \vdots & \vdots & \vdots & \vdots & \ddots & \vdots \\ H_{m11} & \cdots & H_{mn1} & \cdots & \cdots & H_{m1p} & \cdots & H_{mnp} \end{pmatrix}$$

Using the above notation, the operation $u_i = H_{ijk}A_{jk}$ (which is written in indicial notation) can be compactly written as

$$\boldsymbol{u} = \mathrm{mat}_1[\boldsymbol{H}]\mathrm{vec}[\boldsymbol{A}] \tag{43}$$

Let $\boldsymbol{\zeta}$ denote the position vector in the reference element. The row vector containing shape functions is denoted by $\boldsymbol{N}$, which is a function of $\boldsymbol{\zeta}$. The derivatives and hessian of the vector of shape functions with respect to $\boldsymbol{\zeta}$ are, respectively, denoted as $\boldsymbol{DN}$ and $\boldsymbol{DDN}$. That is,

$$(DN)_{ij} = \frac{\partial N_i}{\partial \zeta_j}, \quad \text{and} \quad (DDN)_{ijk} = \frac{\partial^2 N_i}{\partial \zeta_j \partial \zeta_k} \tag{44}$$

The velocity vector is approximated as follows:

$$\boldsymbol{v}(\boldsymbol{x}) = \hat{\boldsymbol{v}}^T \boldsymbol{N}^T(\boldsymbol{\zeta}(\boldsymbol{x})) \tag{45}$$

where $\hat{\boldsymbol{v}}$ denotes the matrix containing nodal velocities. (Note that nodal quantities are denoted with superposed hats.) Based on the isoparametric mapping, the Jacobian matrix and div[$\boldsymbol{J}$] can be written as

$$\boldsymbol{J} = \hat{\boldsymbol{x}}^T(\boldsymbol{DN}) \tag{46}$$

$$\mathrm{div}[\boldsymbol{J}^{-1}] = -\boldsymbol{J}^{-1}\hat{\boldsymbol{x}}^T \mathrm{mat}_1[\boldsymbol{DDN}]\mathrm{vec}[\boldsymbol{J}^{-1}\boldsymbol{J}^{-T}] \tag{47}$$

where $\hat{\boldsymbol{x}}$ is a matrix containing nodal coordinates. For convenience, let us define

$$\boldsymbol{y} = (\boldsymbol{I} - \boldsymbol{DN}\boldsymbol{J}^{-T}\hat{\boldsymbol{x}}^T)\mathrm{mat}_1[\boldsymbol{DDN}]\mathrm{vec}[\boldsymbol{J}^{-1}\boldsymbol{J}^{-T}] \tag{48}$$

Using the above notation, a finite element approximation for $\Delta\boldsymbol{v}$ can be compactly written as

$$\Delta\boldsymbol{v} = (\boldsymbol{y}^T \odot \boldsymbol{I})\mathrm{vec}[\hat{\boldsymbol{v}}^T] \tag{49}$$

Similarly, one can approximate $\Delta\boldsymbol{w}$.

7.2. **Conventional way of deriving a stabilized formulation for transient Brinkman equation.** We define the relevant function spaces, following Hughes [48], as follows (for a more rigorous definition of the function spaces, the reader may consult Reference [49]):

$$\mathcal{V}^t := \left\{\boldsymbol{v}(\cdot,t) \mid \boldsymbol{v}(\cdot,t) \in (H^1(\Omega))^{nd}, \boldsymbol{v}(\boldsymbol{x},t) = \boldsymbol{v}^{\mathrm{p}}(\boldsymbol{x},t) \text{ on } \Gamma\right\} \tag{50a}$$

$$\mathcal{P}^t := \left\{p(\cdot,t) \in L^2(\Omega) \mid \int_\Omega p \; \mathrm{d}\Omega = 0\right\} \tag{50b}$$

$$\mathcal{Q}^t := \left\{p(\cdot,t) \in H^1(\Omega) \mid \int_\Omega p \; \mathrm{d}\Omega = 0\right\} \tag{50c}$$

We again start with the classical mixed formulation for the governing equations: Find $\boldsymbol{v}(\boldsymbol{x},t) \in \mathcal{V}^t$ and $p(\boldsymbol{x},t) \in \mathcal{P}^t$ such that we have

$$\left(\boldsymbol{w}; \rho\frac{\partial \boldsymbol{v}}{\partial t}\right) + (\boldsymbol{w}; \alpha\boldsymbol{v}) + (\mathrm{grad}[\boldsymbol{w}]; \mu\mathrm{grad}[\boldsymbol{v}]) - (\mathrm{div}[\boldsymbol{w}]; p) - (q; \mathrm{div}[\boldsymbol{v}])$$
$$= (\boldsymbol{w}; \rho\boldsymbol{b}) \quad \forall \boldsymbol{w}(\boldsymbol{x}) \in \mathcal{W},\ q(\boldsymbol{x}) \in \mathcal{P} \tag{51}$$

We decompose the velocity as follows:

$$\boldsymbol{v}(\boldsymbol{x},t) = \boldsymbol{v}_c(\boldsymbol{x},t) + \boldsymbol{v}_f(\boldsymbol{x}) \tag{52}$$

As pointed out earlier in Section 1, the above assumption that $\boldsymbol{v}_f(\boldsymbol{x})$ is independent of time, while mathematically plausible, is philosophically undesirable and does not seem to have a physical basis. Substituting (52) into (51), we obtain the coarse-scale and fine-scale subproblems, which can be written as follows: Find $\boldsymbol{v}_c(\boldsymbol{x},t) \in \mathcal{V}^t$, $\boldsymbol{v}_f(\boldsymbol{x}) \in \mathcal{V}_f$ and $p(\boldsymbol{x},t) \in \mathcal{P}^t$ such that

$$\begin{aligned}
&\left(\boldsymbol{w}_c; \rho\frac{\partial \boldsymbol{v}_c}{\partial t}\right) + (\boldsymbol{w}_c; \alpha\boldsymbol{v}_c) + (\boldsymbol{w}_c; \alpha\boldsymbol{v}_f) + (\mathrm{grad}[\boldsymbol{w}_c]; \mu\ \mathrm{grad}[\boldsymbol{v}_c]) + (\mathrm{grad}[\boldsymbol{w}_c]; \mu\ \mathrm{grad}[\boldsymbol{v}_f]) \\
&- (\mathrm{div}[\boldsymbol{w}_c]; p) - (q; \mathrm{div}[\boldsymbol{v}_c]) - (q; \mathrm{div}[\boldsymbol{v}_f]) = (\boldsymbol{w}_c; \rho\boldsymbol{b}) \quad \forall \boldsymbol{w}_c(\boldsymbol{x}) \in \mathcal{W}_c,\ q(\boldsymbol{x}) \in \mathcal{P} \qquad (53a) \\
&\left(\boldsymbol{w}_f; \rho\frac{\partial \boldsymbol{v}_c}{\partial t}\right) + (\boldsymbol{w}_f; \alpha\boldsymbol{v}_c) + (\boldsymbol{w}_f; \alpha\boldsymbol{v}_f) + (\mathrm{grad}[\boldsymbol{w}_f]; \mu\ \mathrm{grad}[\boldsymbol{v}_c]) + (\mathrm{grad}[\boldsymbol{w}_f]; \mu\ \mathrm{grad}[\boldsymbol{v}_f]) \\
&- (\mathrm{div}[\boldsymbol{w}_f]; p) = (\boldsymbol{w}_f; \rho\boldsymbol{b}) \quad \forall \boldsymbol{w}_f(\boldsymbol{x}) \in \mathcal{W}_f \qquad (53b)
\end{aligned}$$

The remaining part of the derivation is as follows: we interpolate the fine-scale variables ($\boldsymbol{v}_f(\boldsymbol{x})$ and $\boldsymbol{w}_f(\boldsymbol{x})$) in terms of the bubble function $b^e(\boldsymbol{x})$. We then solve the fine-scale problem to obtain the fine-scale component of the velocity in terms of the coarse-scale component of the velocity and pressure. The obtained expression for the fine-scale component of the velocity is then substituted in the coarse-scale subproblem (53a) to eliminate the fine-scale degrees-of-freedom and obtain a stabilized mixed formulation. For convenience, dropping the suffix $c$, the final form of the stabilized mixed formulation reads: Find $\boldsymbol{v}(\boldsymbol{x},t) \in \mathcal{V}^t$ and $p(\boldsymbol{x},t) \in \mathcal{Q}^t$ such that we have

$$\begin{aligned}
&\left(\boldsymbol{w}; \rho\frac{\partial \boldsymbol{v}}{\partial t}\right) + (\boldsymbol{w}; \alpha\boldsymbol{v}) + (\mathrm{grad}[\boldsymbol{w}]; \mu\ \mathrm{grad}[\boldsymbol{v}]) - (\mathrm{div}[\boldsymbol{w}]; p) - (q; \mathrm{div}[\boldsymbol{v}]) \\
&- \sum_{e=1}^{Nele} \left(\alpha\boldsymbol{w} + \mathrm{grad}[q] - \mu\mathrm{div}[\mathrm{grad}[\boldsymbol{w}]; \tilde{\boldsymbol{\tau}}(\boldsymbol{x})\rho\frac{\partial \boldsymbol{v}}{\partial t}\right)_{\Omega^e} \\
&- \sum_{e=1}^{Nele} (\alpha\boldsymbol{w} + \mathrm{grad}[q] - \mu\mathrm{div}[\mathrm{grad}[\boldsymbol{w}]; \tilde{\boldsymbol{\tau}}(\boldsymbol{x})\,(\alpha\boldsymbol{v} + \mathrm{grad}[p] - \mu\mathrm{div}[\mathrm{grad}[\boldsymbol{v}]))_{\Omega^e} \\
&= (\boldsymbol{w}; \rho\boldsymbol{b}) - \sum_{e=1}^{Nele} (\alpha\boldsymbol{w} + \mathrm{grad}[q] - \mu\mathrm{div}[\mathrm{grad}[\boldsymbol{w}]; \tilde{\boldsymbol{\tau}}(\boldsymbol{x})\rho\boldsymbol{b})_{\Omega^e} \quad \forall \boldsymbol{w}(\boldsymbol{x}) \in \mathcal{W},\ q(\boldsymbol{x}) \in \mathcal{Q} \qquad (54)
\end{aligned}$$

where the stabilization parameter in this case is defined as follows:

$$\tilde{\boldsymbol{\tau}}(\boldsymbol{x}) = b^e(\boldsymbol{x}) \left[\int_{\Omega^e} b^e(\boldsymbol{x})\ \mathrm{d}\Omega\right] \left[\int_{\Omega^e} \left(\mu\ \|\mathrm{grad}[b^e]\|^2 + \alpha(b^e)^2\right)\ \mathrm{d}\Omega\right]^{-1} \tag{55}$$

A semi-discretization of equation (54) yields a system of ordinary differential equations (ODEs) of the following general form [48]:

$$\begin{bmatrix} \mathbf{M_{vv}} & \mathbf{M_{vp}} \\ \mathbf{M_{pv}} & \mathbf{0} \end{bmatrix} \begin{bmatrix} \dot{\mathbf{v}} \\ \dot{\mathbf{p}} \end{bmatrix} + \begin{bmatrix} \mathbf{K_{vv}} & \mathbf{K_{vp}} \\ \mathbf{K_{pv}} & \mathbf{K_{pp}} \end{bmatrix} \begin{bmatrix} \mathbf{v} \\ \mathbf{p} \end{bmatrix} = \begin{bmatrix} \mathbf{f_v} \\ \mathbf{f_p} \end{bmatrix} \tag{56}$$

The above system of ODEs can be solved by employing any of the known time stepping schemes (e.g., Runge-Kutta methods, generalized trapezoidal family of time stepping schemes).

## ACKNOWLEDGMENTS

The research reported herein was supported by Texas Engineering Experiment Station (TEES). This support is gratefully acknowledged. The opinions expressed in this paper are those of the authors and do not necessarily reflect that of the sponsor.

Shriram Srinivasan, Graduate student, Department of Mechanical Engineering, Texas A&M University, College Station, Texas 77843.

*E-mail address*: `shrirams@tamu.edu`

Correspondence to: Dr. Kalyana Babu Nakshatrala, Department of Mechanical Engineering, 216 Engineering/Physics Building, Texas A&M University, College Station, Texas 77843. TEL:+1-979-845-1292

*E-mail address*: `knakshatrala@tamu.edu`

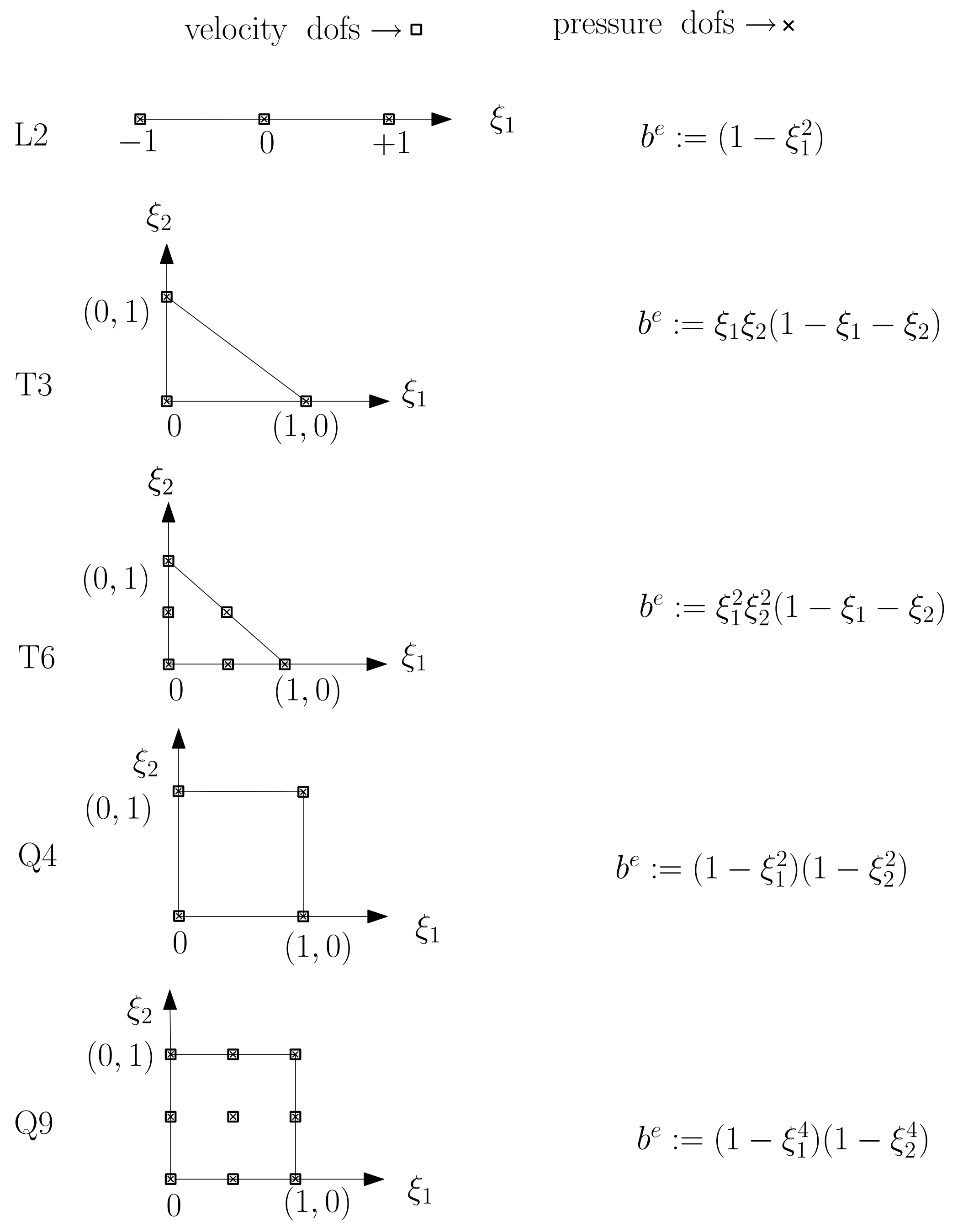

Figure 2. This figure shows the bubble functions for various finite elements. The coordinate axes of the reference/master element are denoted by $\xi_1$ and $\xi_2$. The locations of the degrees-of-freedom for velocity and pressure nodes in each element are also indicated.

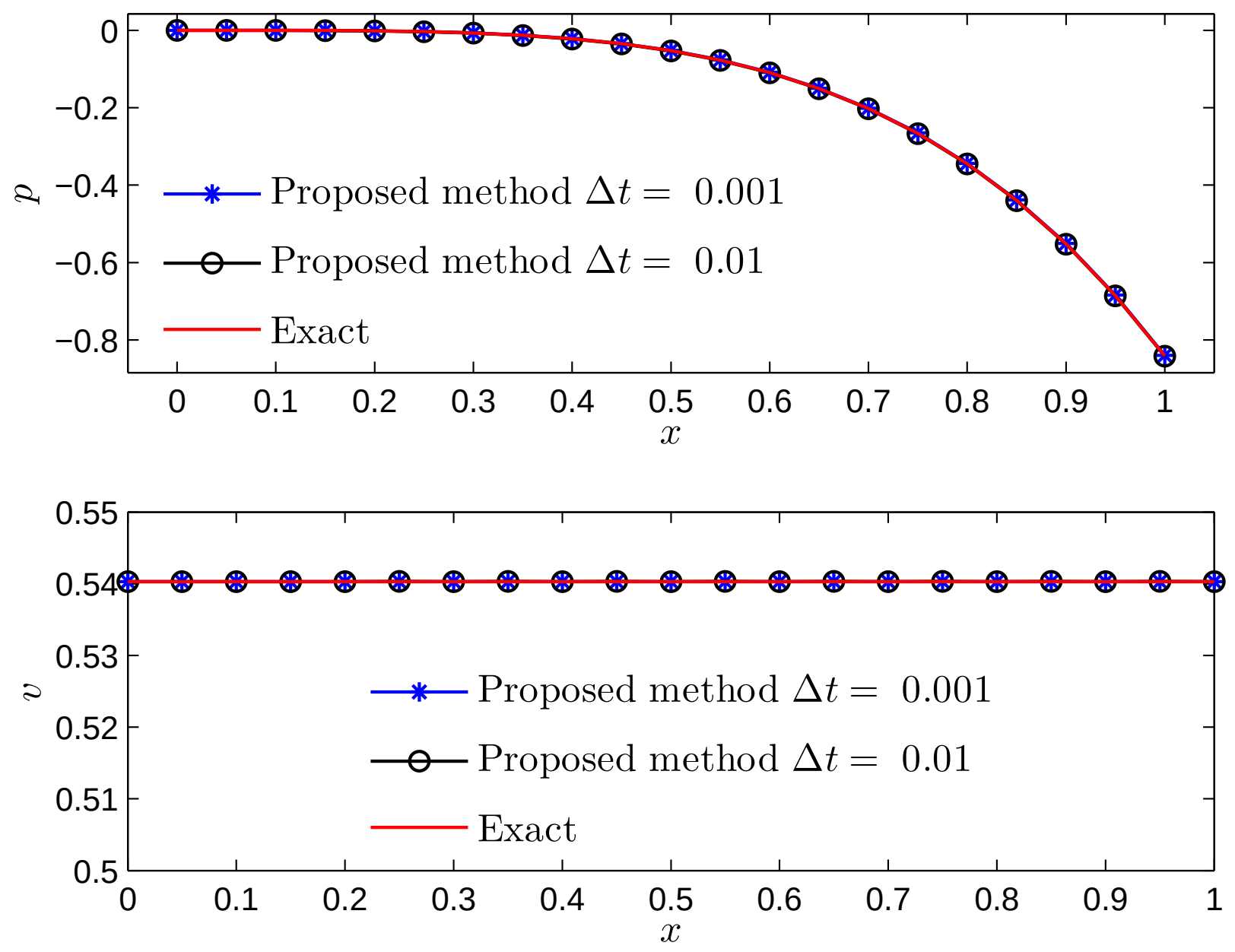

FIGURE 3. One-dimensional test problem (unsteady Darcy equation): In this figure the obtained numerical solution is compared with the analytical solution for the time instant at $T = 1$. The computational mesh consisted of 20 equal-sized two-node linear finite elements.

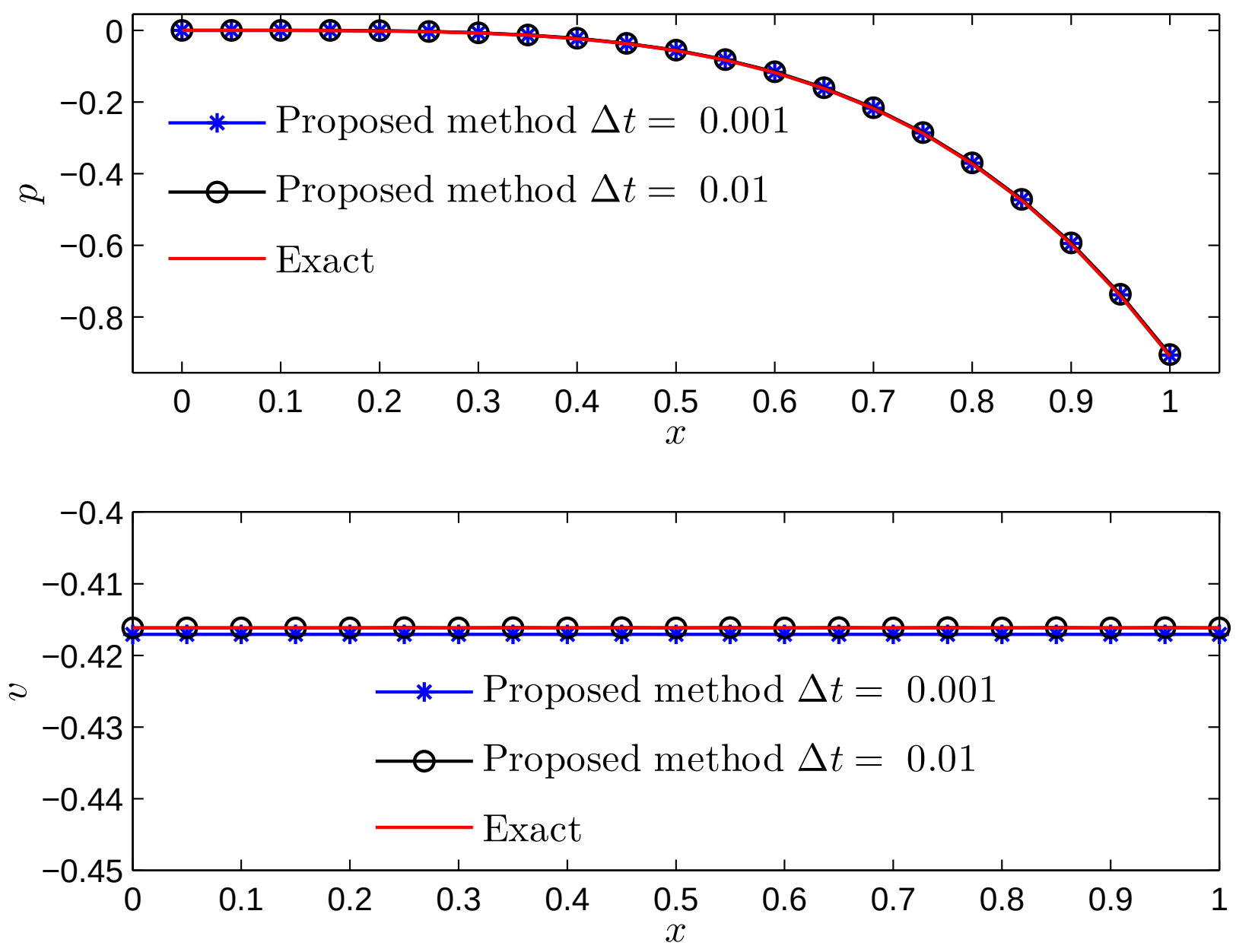

FIGURE 4. One-dimensional test problem (unsteady Darcy equation): In this figure the obtained numerical solution is compared with the analytical solution for the time instant at $T = 2$. The computational mesh consisted of 20 equal-sized two-node linear finite elements.

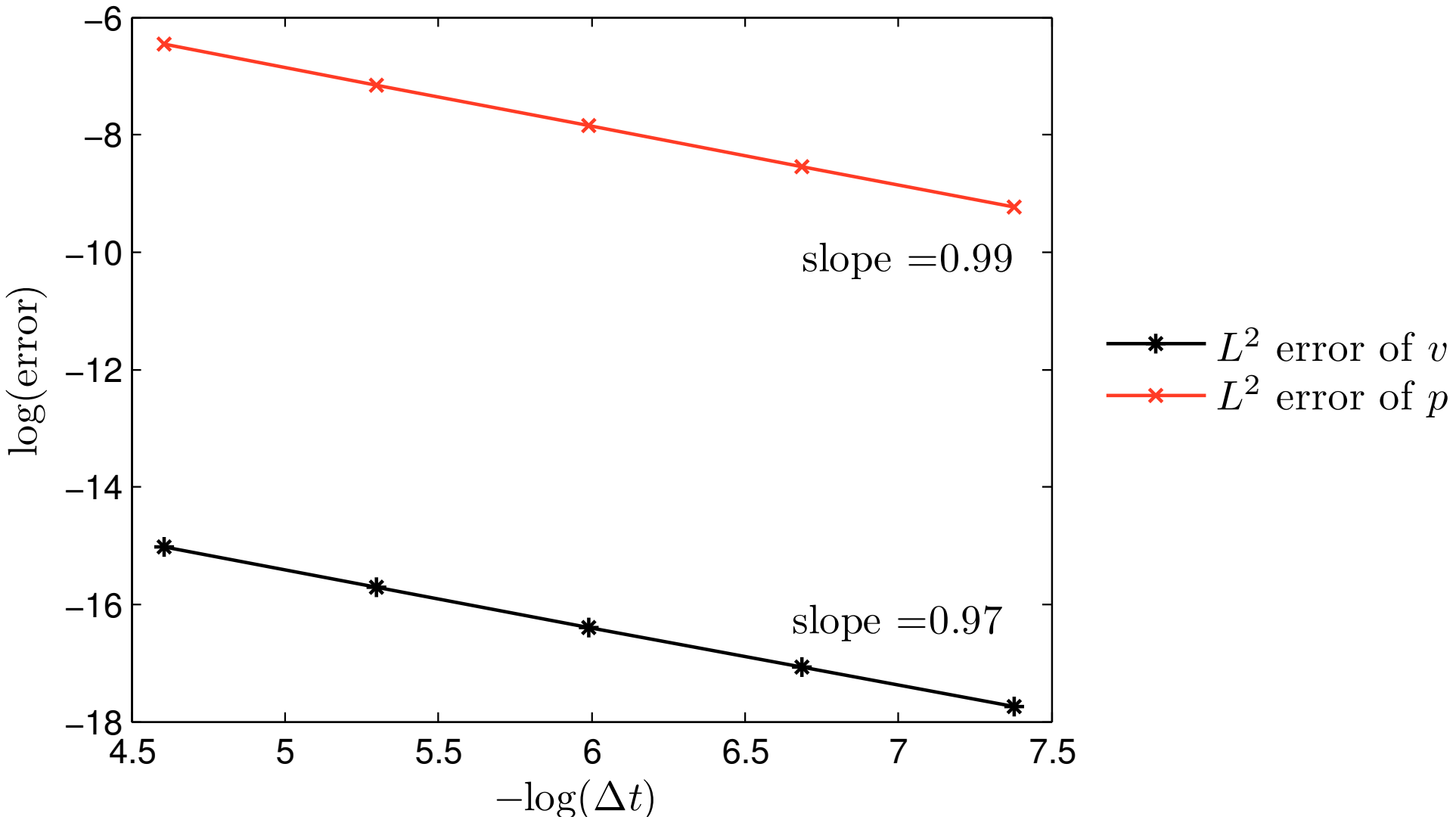

(a) Temporal convergence using 128 elements after time $T = 1$

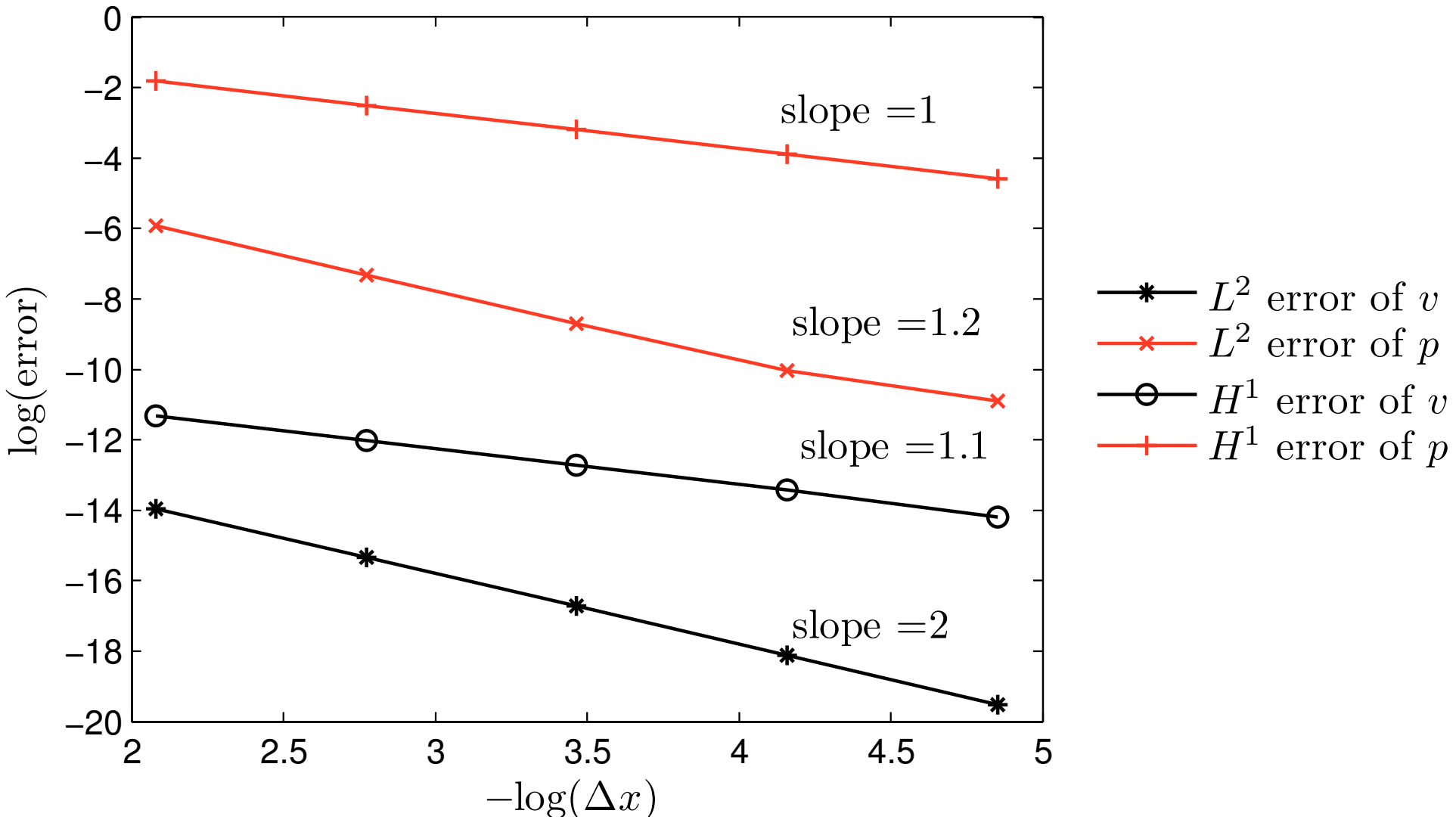

(b) Spatial convergence using $\Delta t = 10^{-4}$ after time $T = 1$

FIGURE 5. One-dimensional test problem (unsteady Darcy equation): Temporal and spatial convergence. Note that the rate of convergence is in good agreement with the theory.

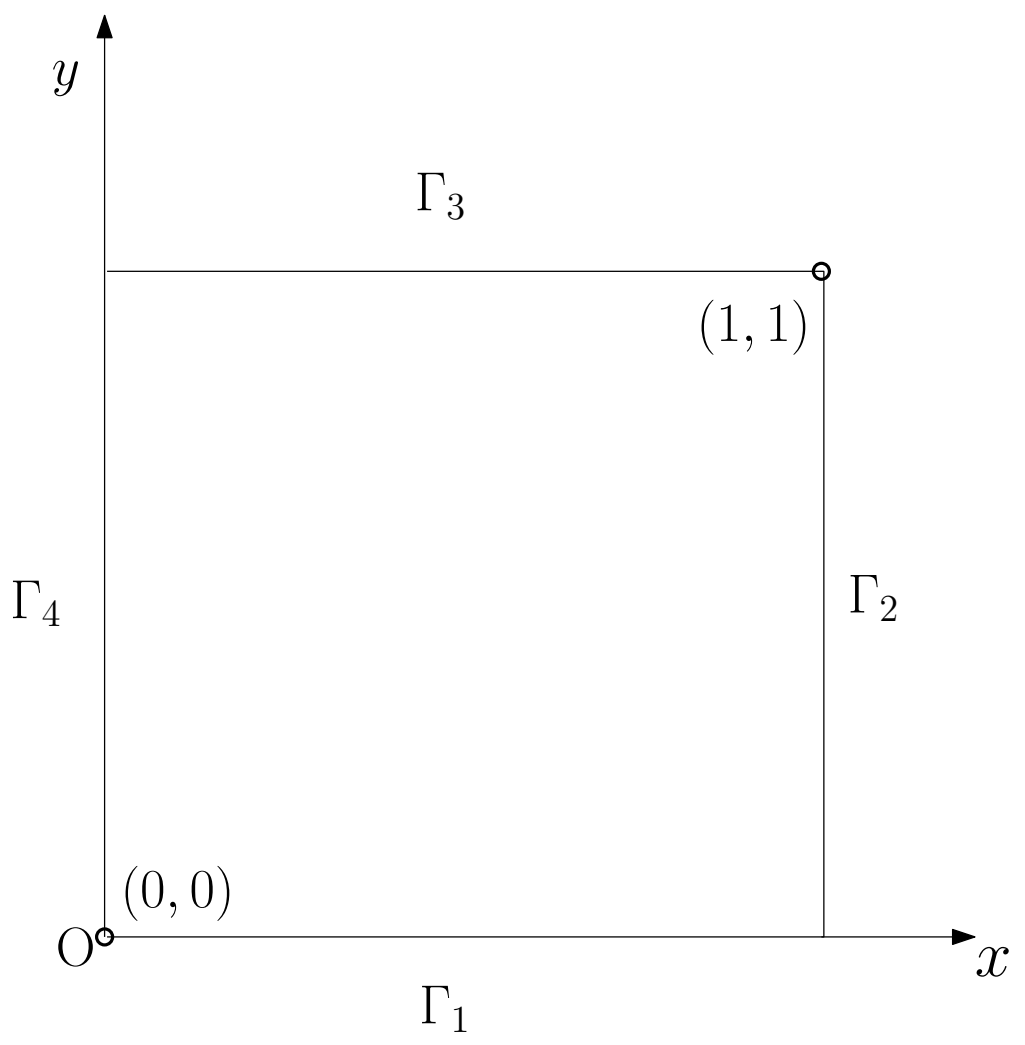

Figure 6. Two-dimensional problem : A pictorial description of the test problem. The computational domain is a bi-unit square. For unsteady Darcy equation: Homogeneous normal components of the velocity are prescribed on the whole boundary (that is, $\psi(\boldsymbol{x}, t) = 0$). The components of the velocity vector are, respectively, denoted by $v_x$ and $v_y$ so that $v_y = 0$ on $\Gamma_1,\ \Gamma_3$ while $v_x = 0$ on $\Gamma_2,\ \Gamma_4$. For unsteady Brinkman equation: Velocity is prescribed completely on the boundary, so that, $v_x = v_y = 0$ on $\Gamma_1,\ \Gamma_3$ while $v_x = 0,\ v_y = (\exp(t) - 1)\,(\cos(2\pi y) - 1)$ on $\Gamma_2,\ \Gamma_4$. To ensure uniqueness of the solution, we enforce pressure to be zero at the origin O in both cases.

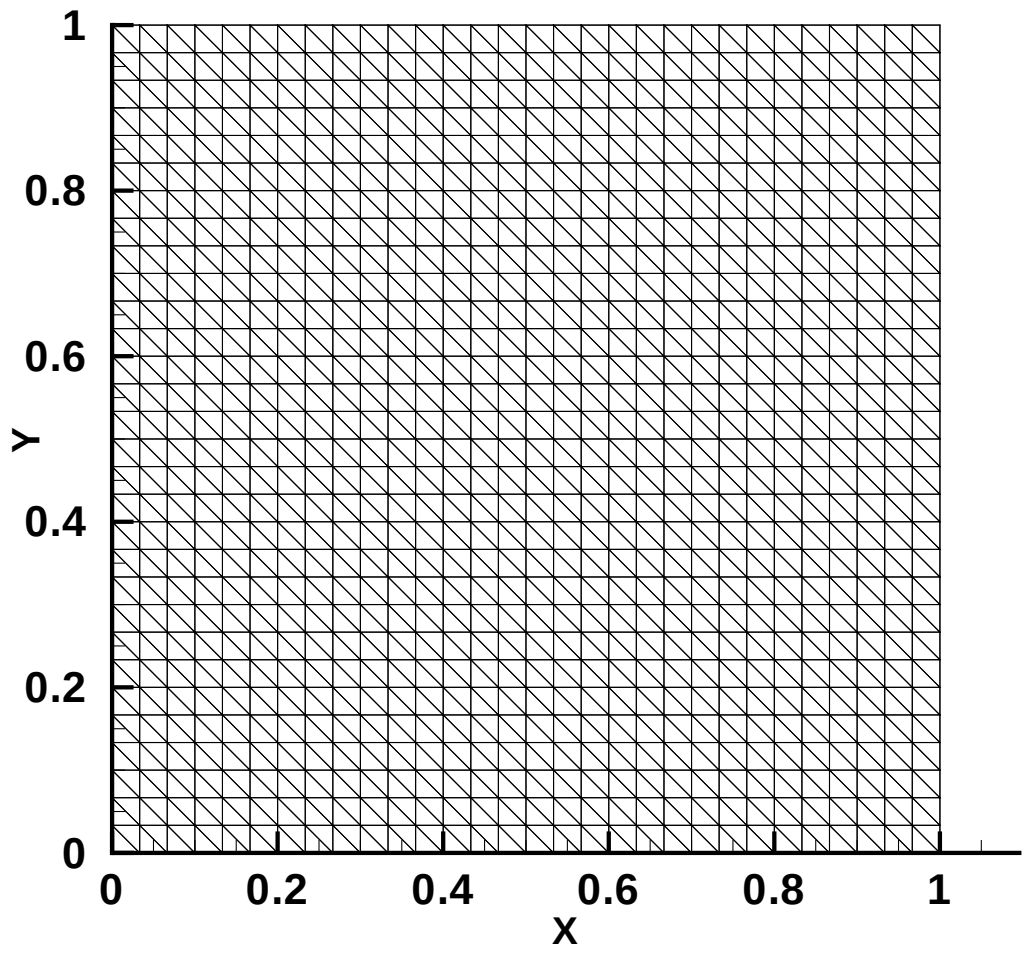

(a) Three-node triangular mesh with 1800 elements

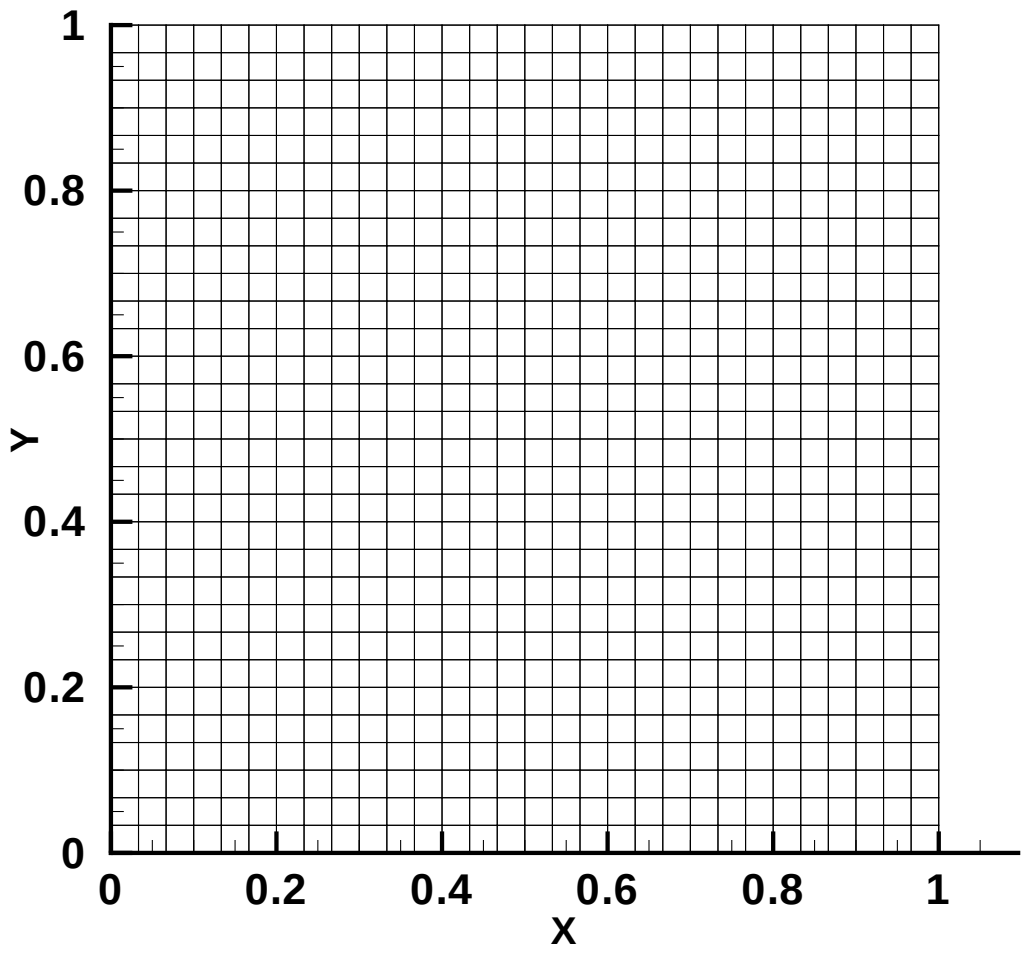

(b) Four-node quadrilateral mesh with 900 elements

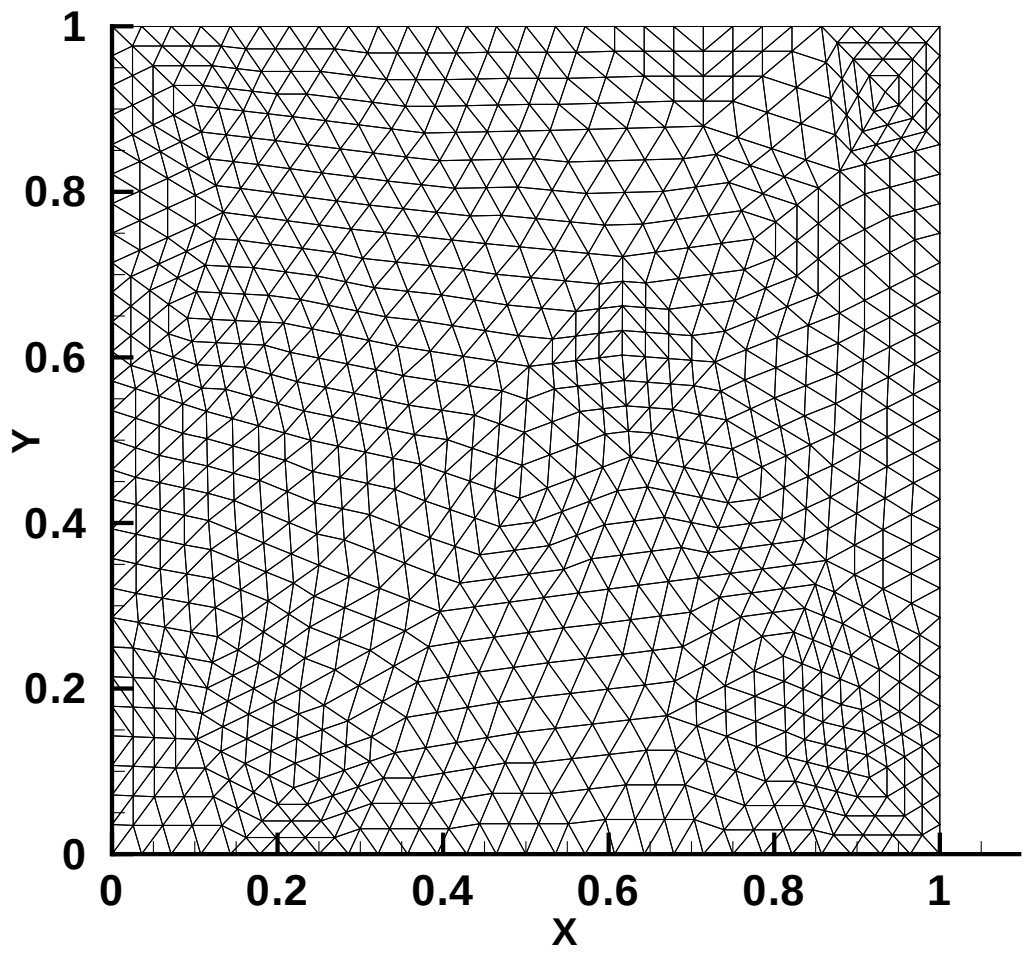

(c) Unstructured three-node triangular mesh with 1792 elements

FIGURE 7. Two-dimensional test problem: Computational meshes used in the numerical simulations.

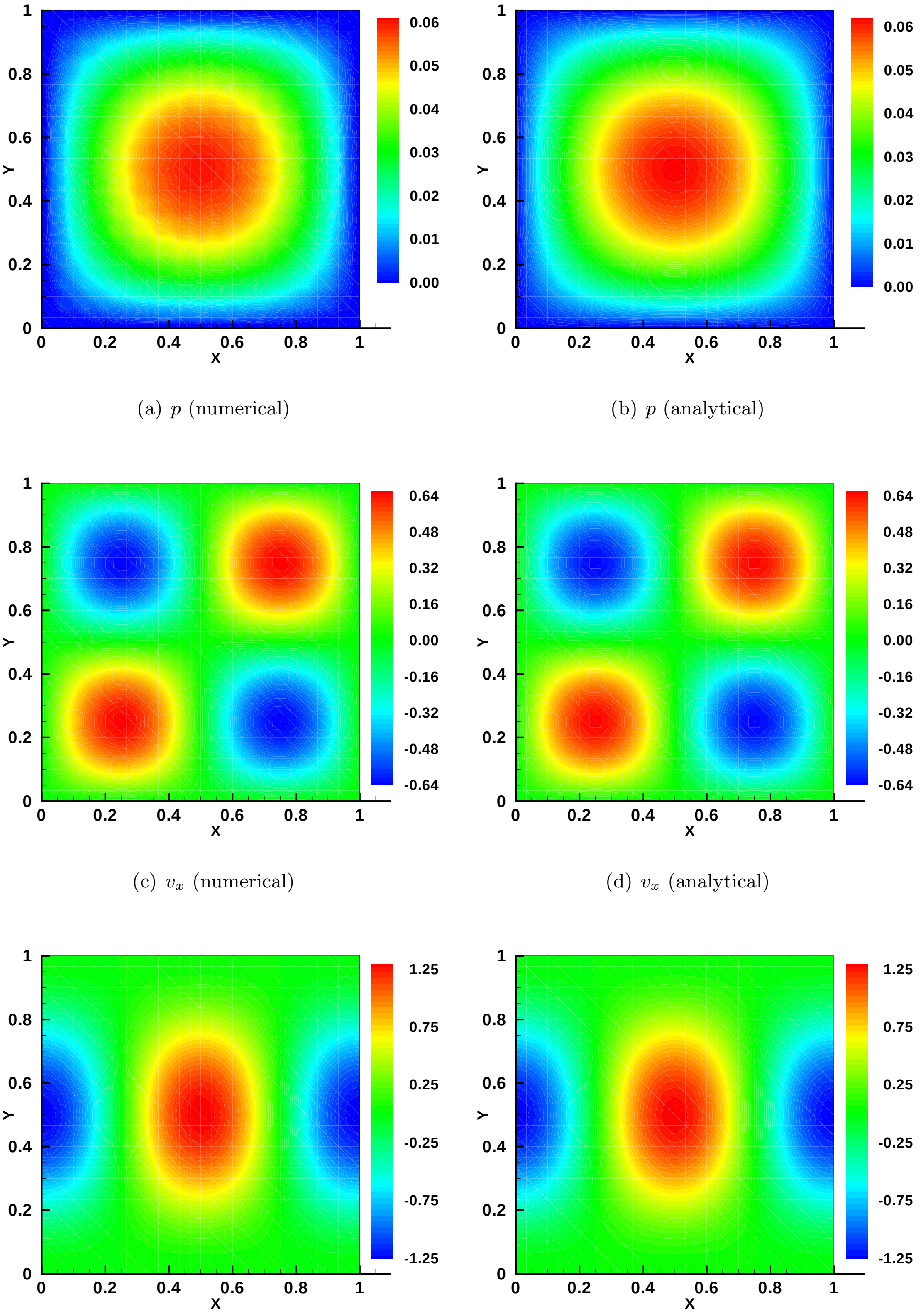

(a) $p$ (numerical) (b) $p$ (analytical)

(c) $v_x$ (numerical) (d) $v_x$ (analytical)

(e) $v_y$ (numerical) (f) $v_y$ (analytical)

FIGURE 8. Two-dimensional test problem (unsteady Darcy equation) with *three-node triangular finite elements*: Comparison of analytical and numerical solutions at final time $T = 0.5$. An equally spaced $31 \times 31$ grid was used with $\Delta t = 10^{-3}$.

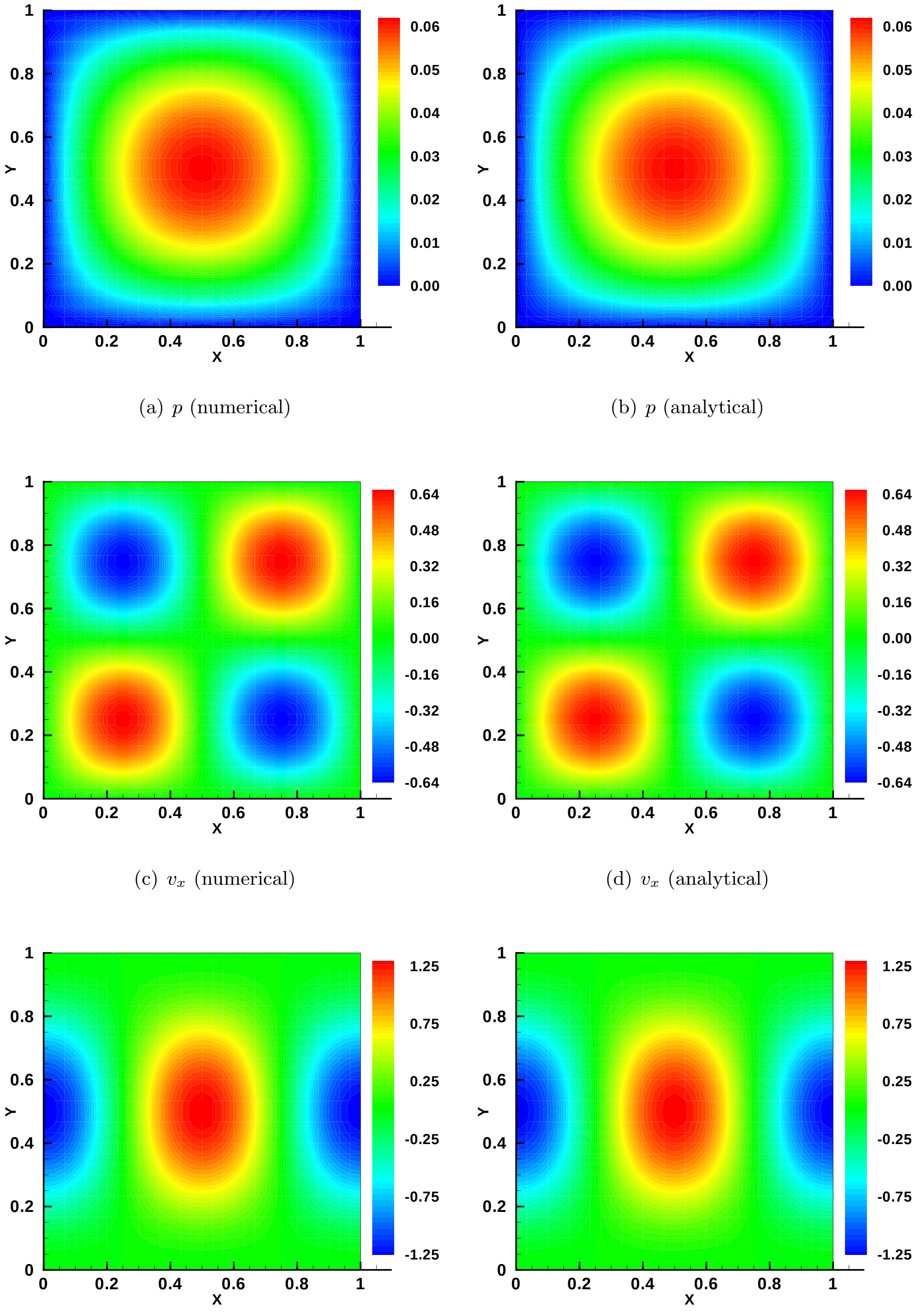

(a) $p$ (numerical)

(b) $p$ (analytical)

(c) $v_x$ (numerical)

(d) $v_x$ (analytical)

(e) $v_y$ (numerical)

(f) $v_y$ (analytical)

FIGURE 9. Two-dimensional test problem (unsteady Darcy equation) with *four-node quadrilateral finite elements*: Comparison of analytical and numerical solutions at final time $T = 0.5$. An equally spaced $31 \times 31$ grid was used with $\Delta t = 10^{-3}$.

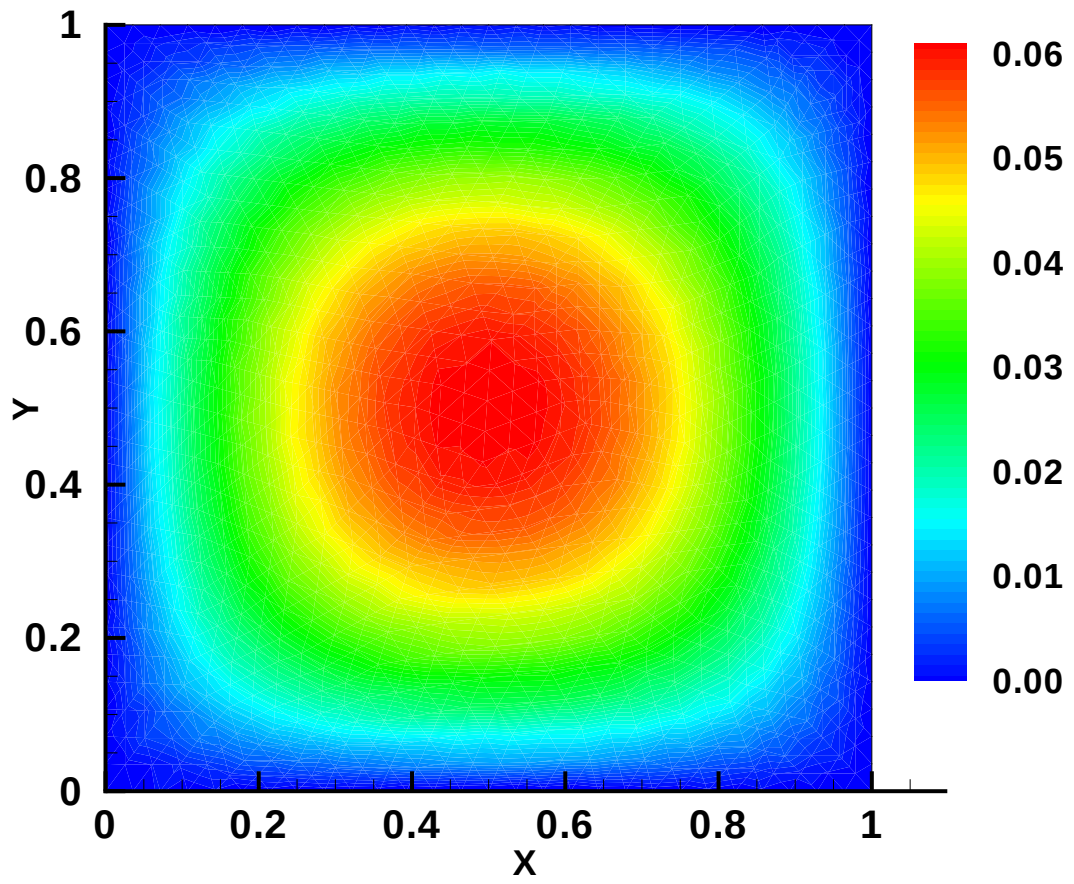

(a) $p$ (numerical)

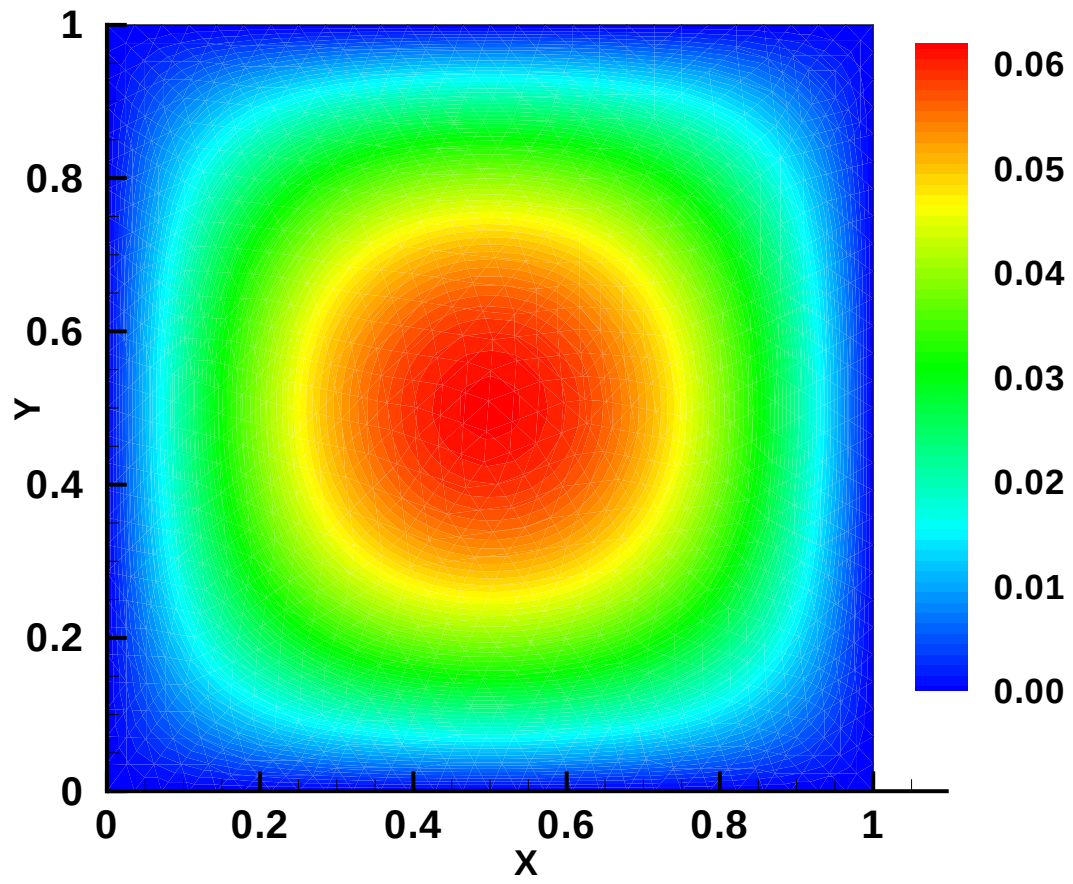

(b) $p$ (analytical)

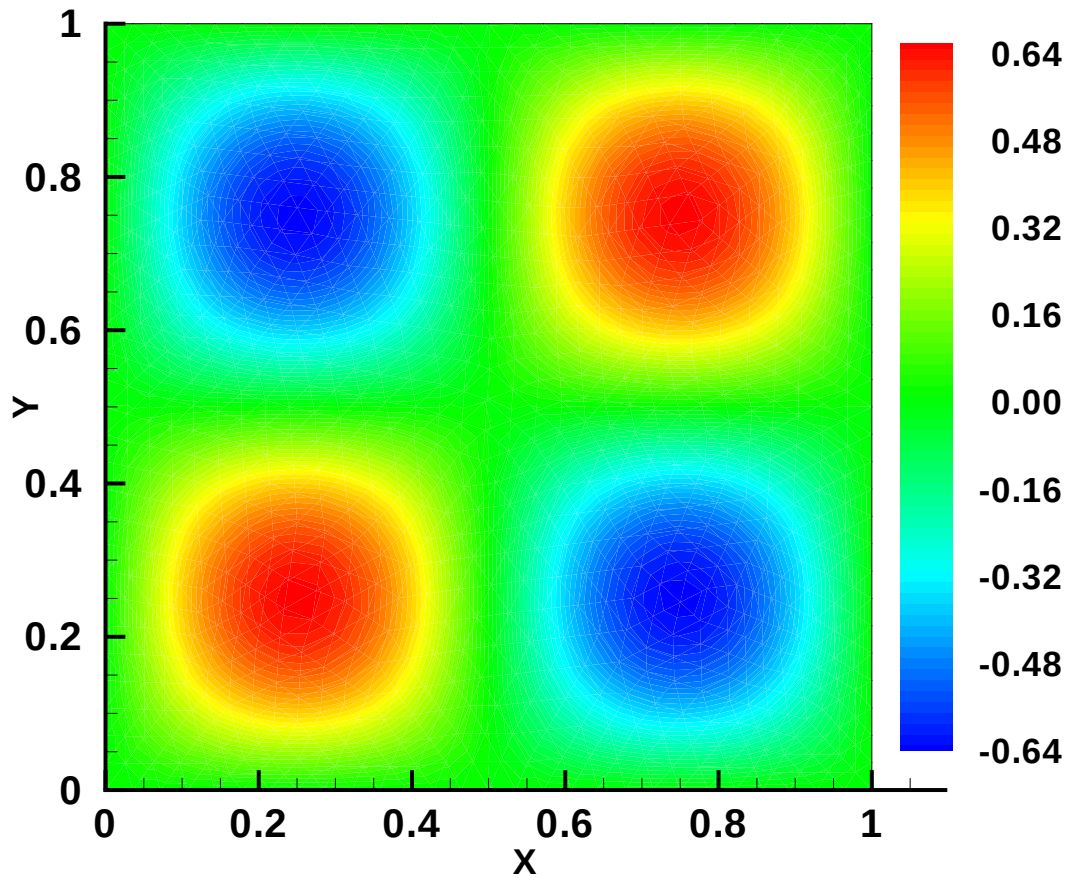

(c) $v_x$ (numerical)

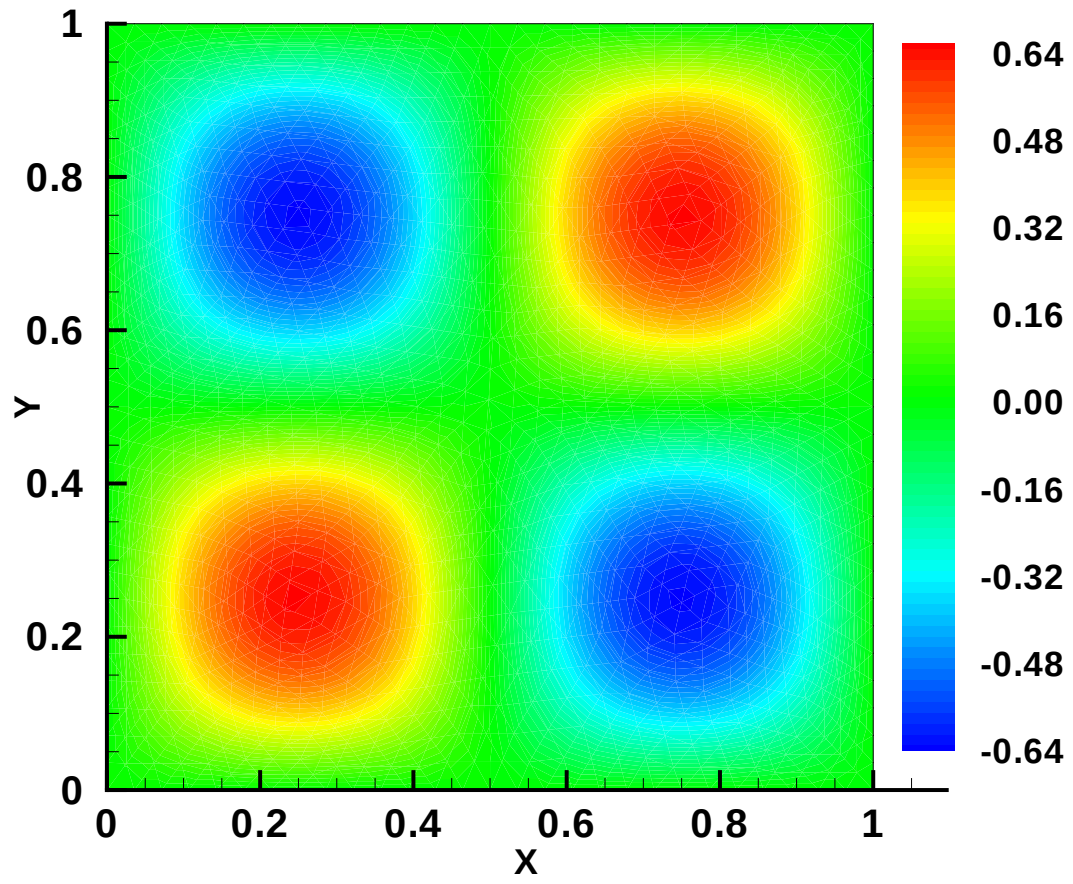

(d) $v_x$ (analytical)

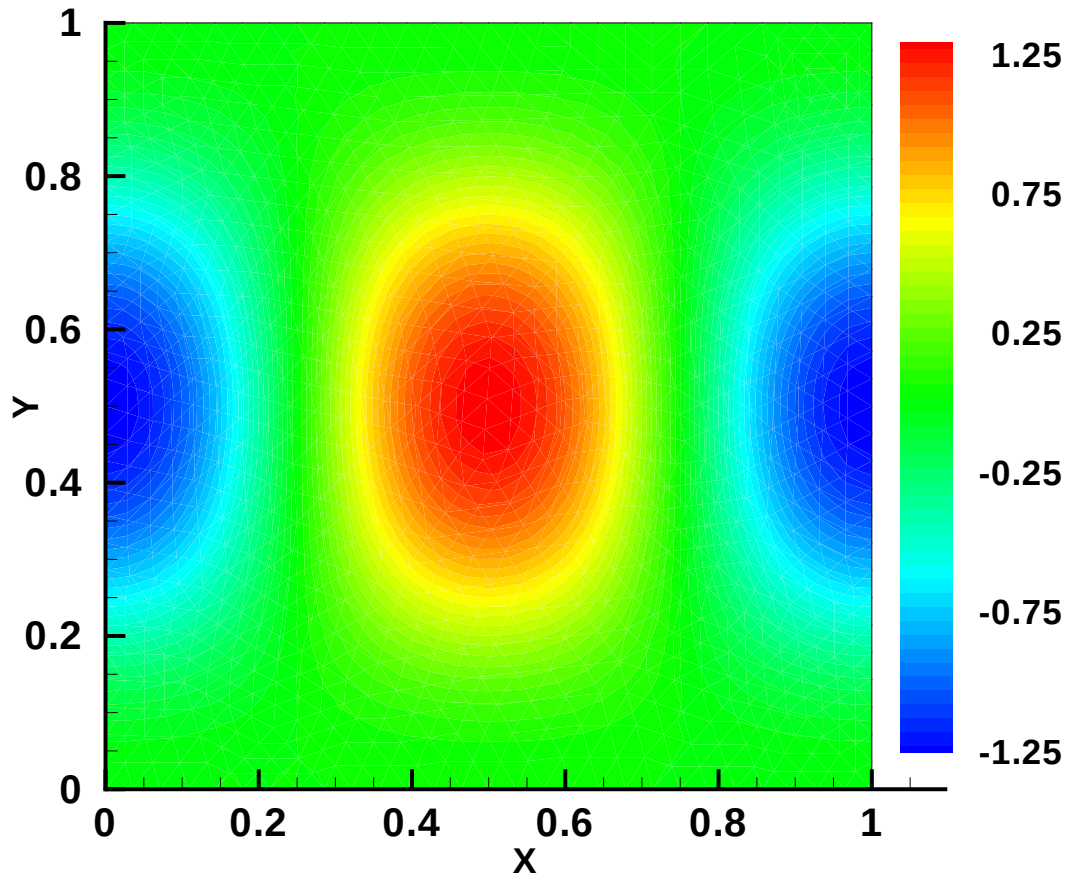

(e) $v_y$ (numerical)

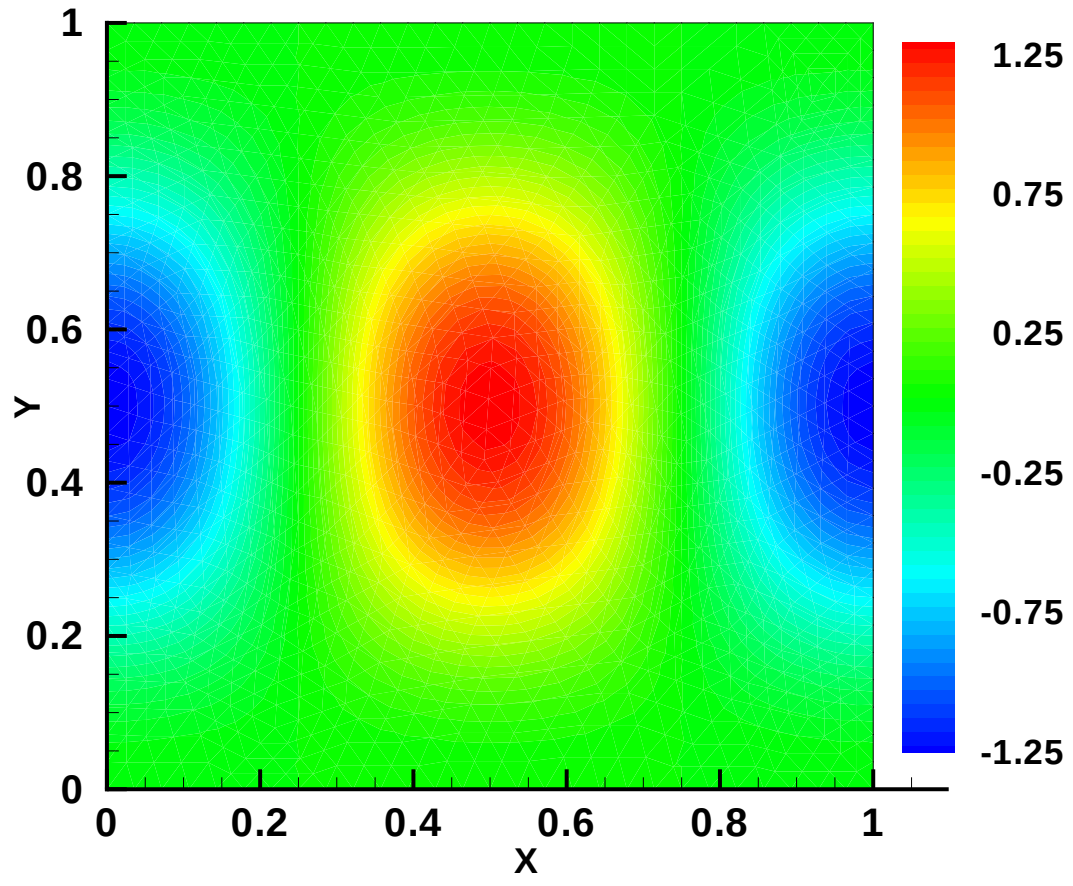

(f) $v_y$ (analytical)

FIGURE 10. Two-dimensional test problem (unsteady Darcy equation) with *unstructured grid of three-node finite elements*: Comparison of analytical and numerical solutions at final time $T = 0.5$. 1792 elements were used with $\Delta t = 10^{-3}$.

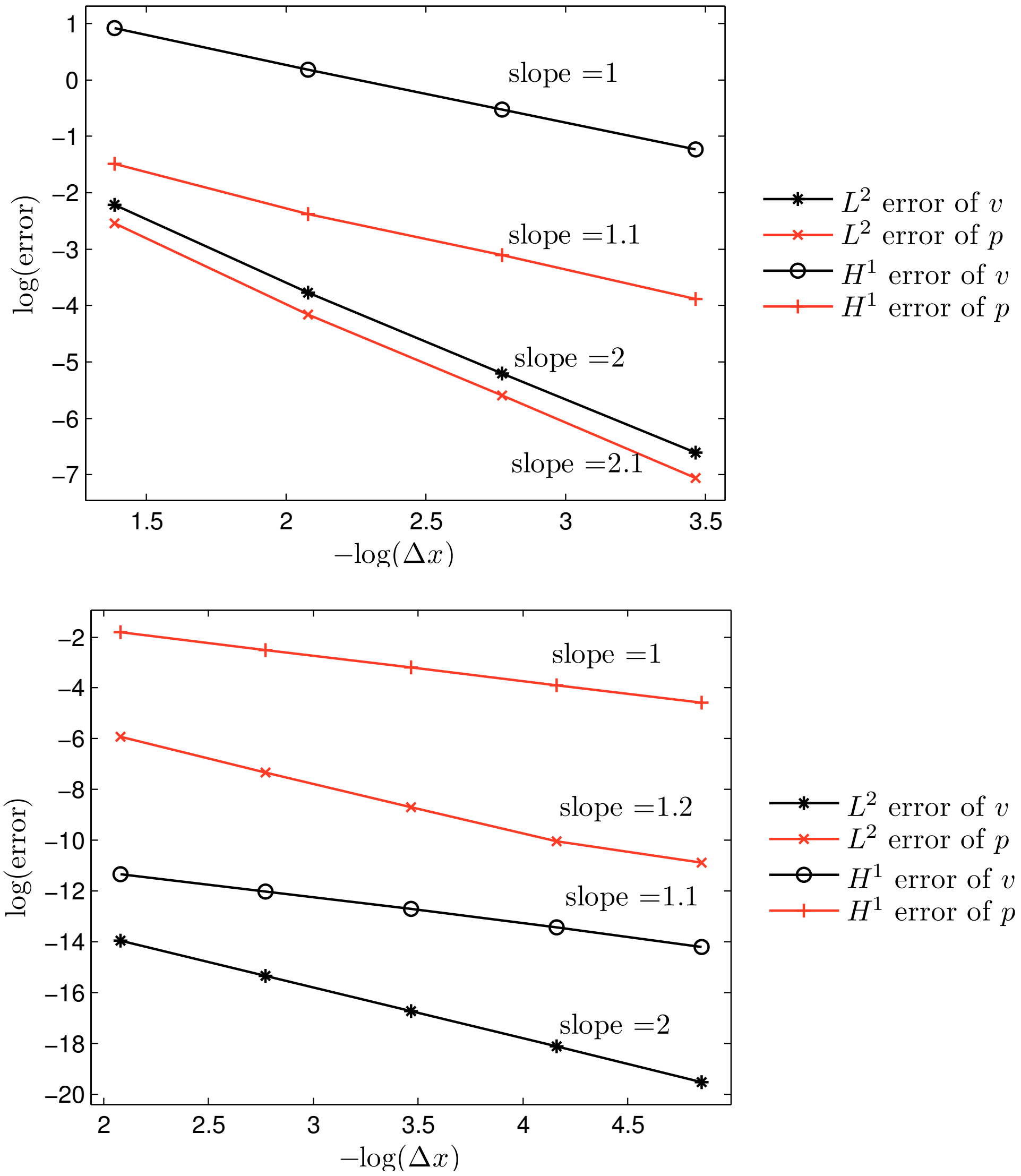

FIGURE 11. Two-dimensional test problem (unsteady Darcy equation): The figure presents spatial convergence of the proposed formulation using four-node quadrilateral elements (top) and three-node triangular elements (bottom). The time step is taken to be $\Delta t = 10^{-3}$, and the total time is taken to be $T = 0.5$.

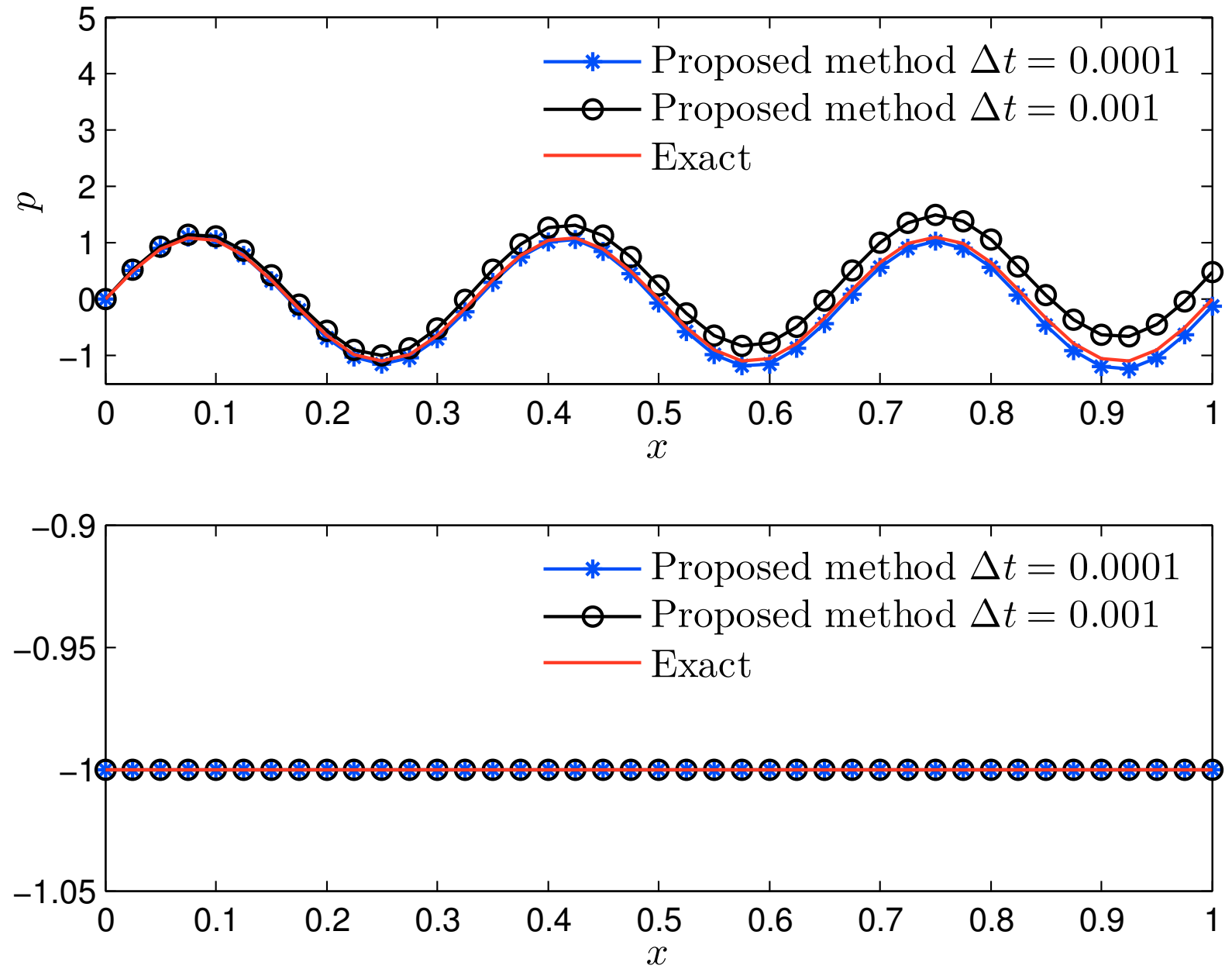

FIGURE 12. One-dimensional test problem (unsteady Brinkman equation): In this figure the obtained numerical solution is compared with the analytical solution for the time instant at $T = 0.1$. The computational mesh consisted of 40 equal-sized two-node linear finite elements.

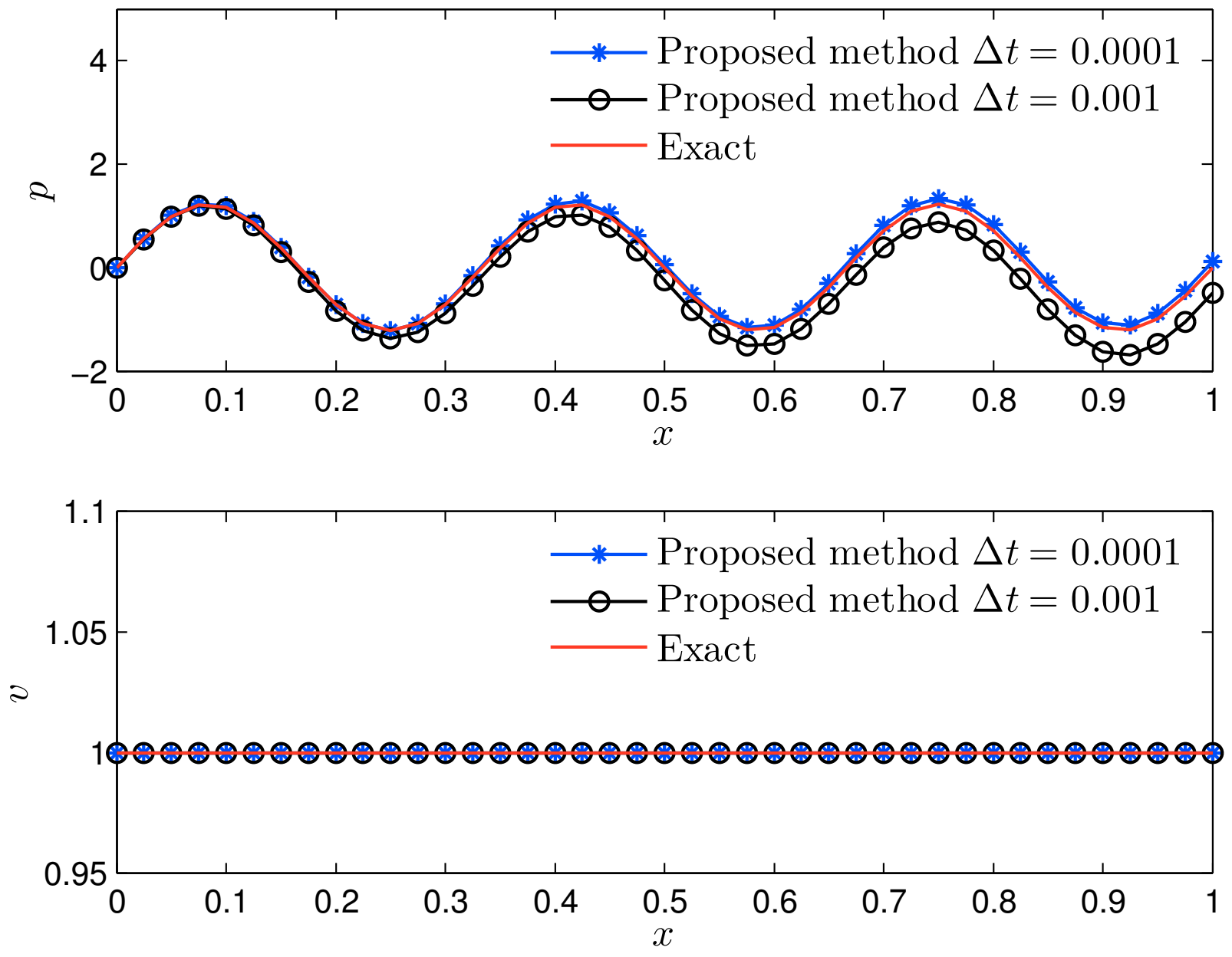

FIGURE 13. One-dimensional test problem (unsteady Brinkman equation): In this figure the obtained numerical solution is compared with the analytical solution for the time instant at $T = 0.2$. The computational mesh consisted of 40 equal-sized two-node linear finite elements.

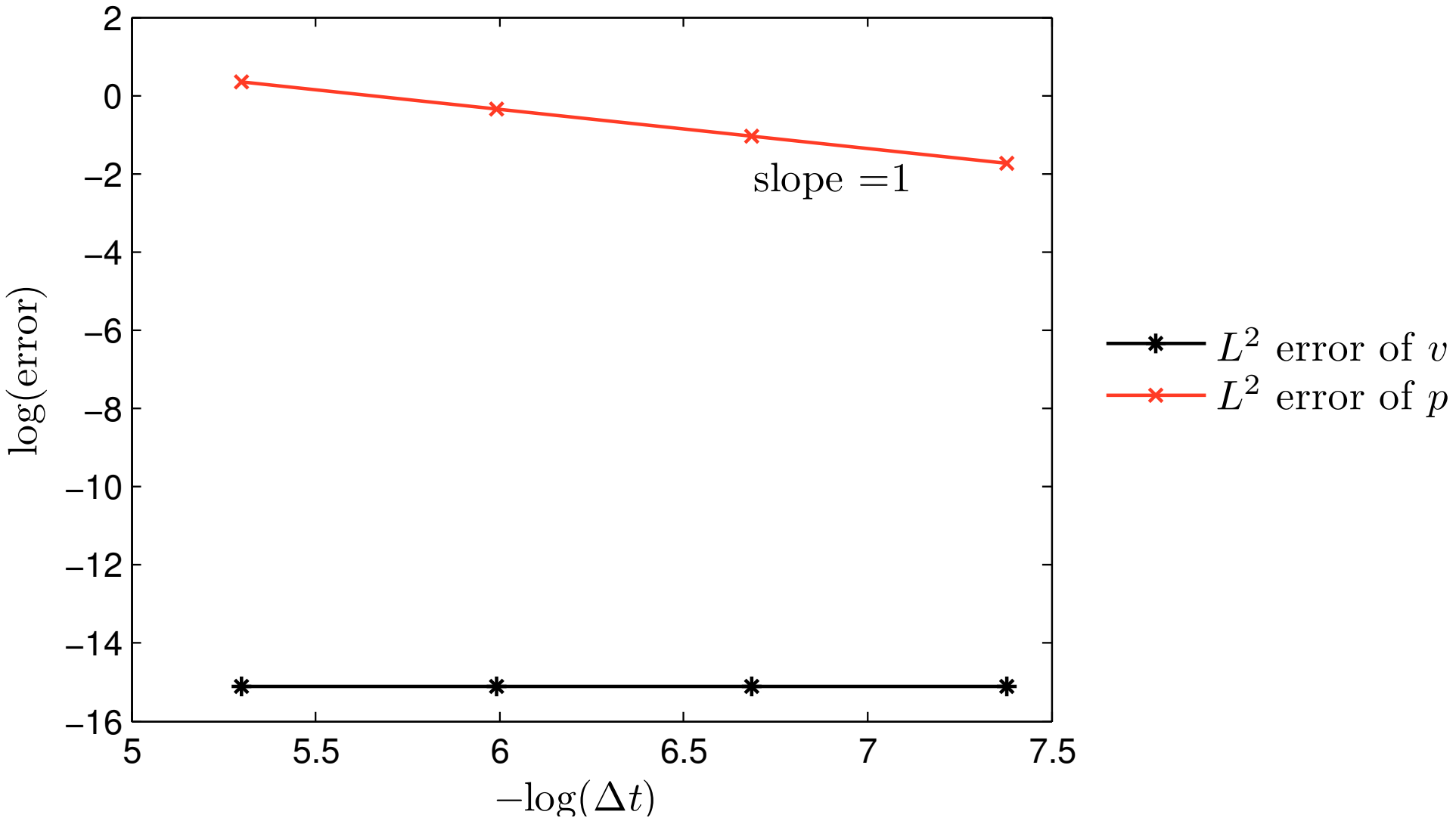

(a) Temporal convergence using 128 elements after time $T = 1$

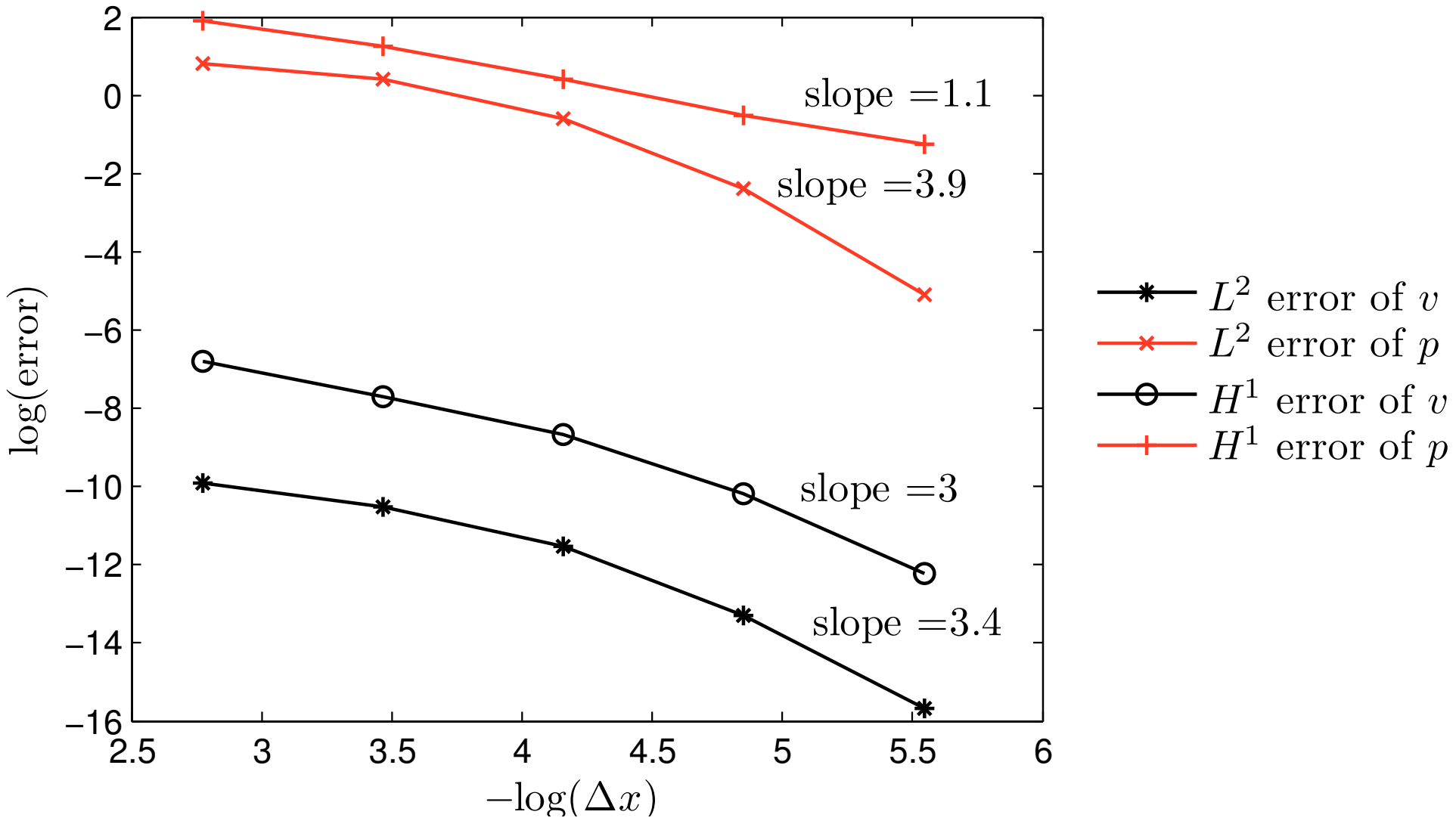

(b) Spatial convergence using $\Delta t = 10^{-5}$ after time $T = 0.01$

FIGURE 14. One-dimensional test problem (unsteady Brinkman equation): Temporal and spatial convergence. Note that the rate of convergence is in good agreement with the theory.

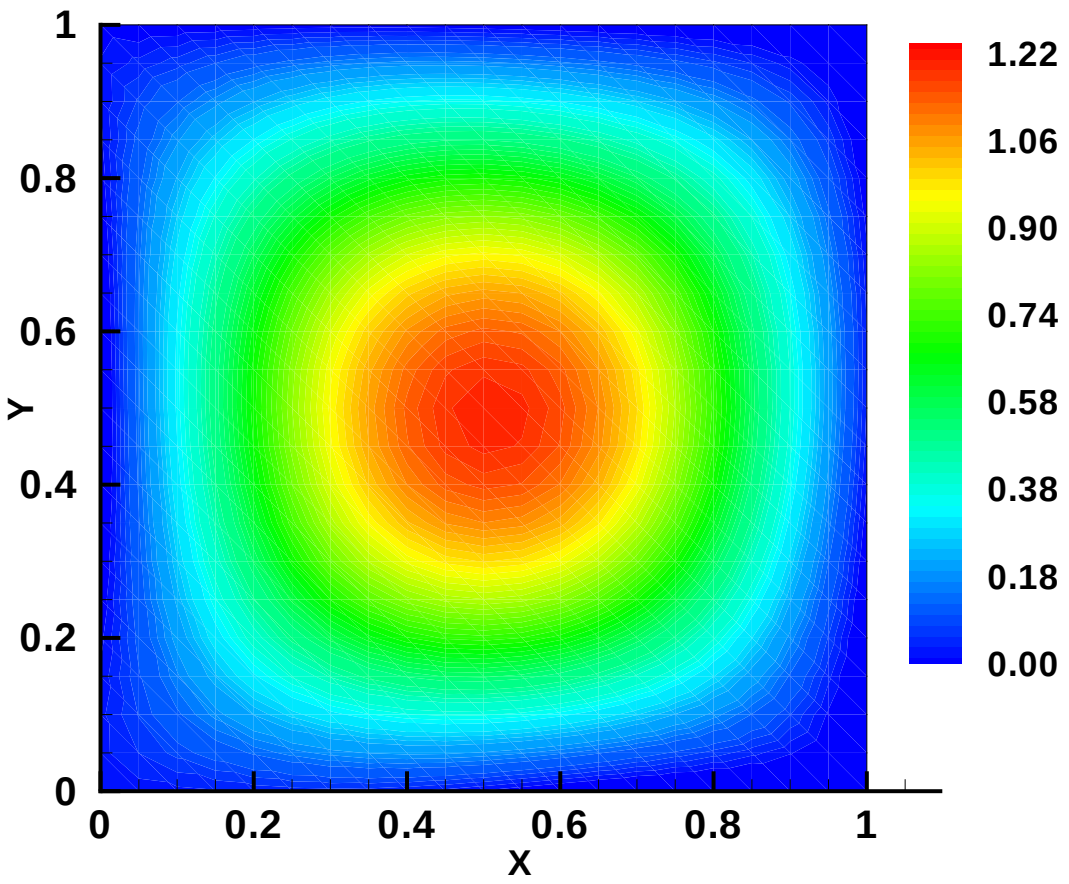

(a) $p$ (numerical)

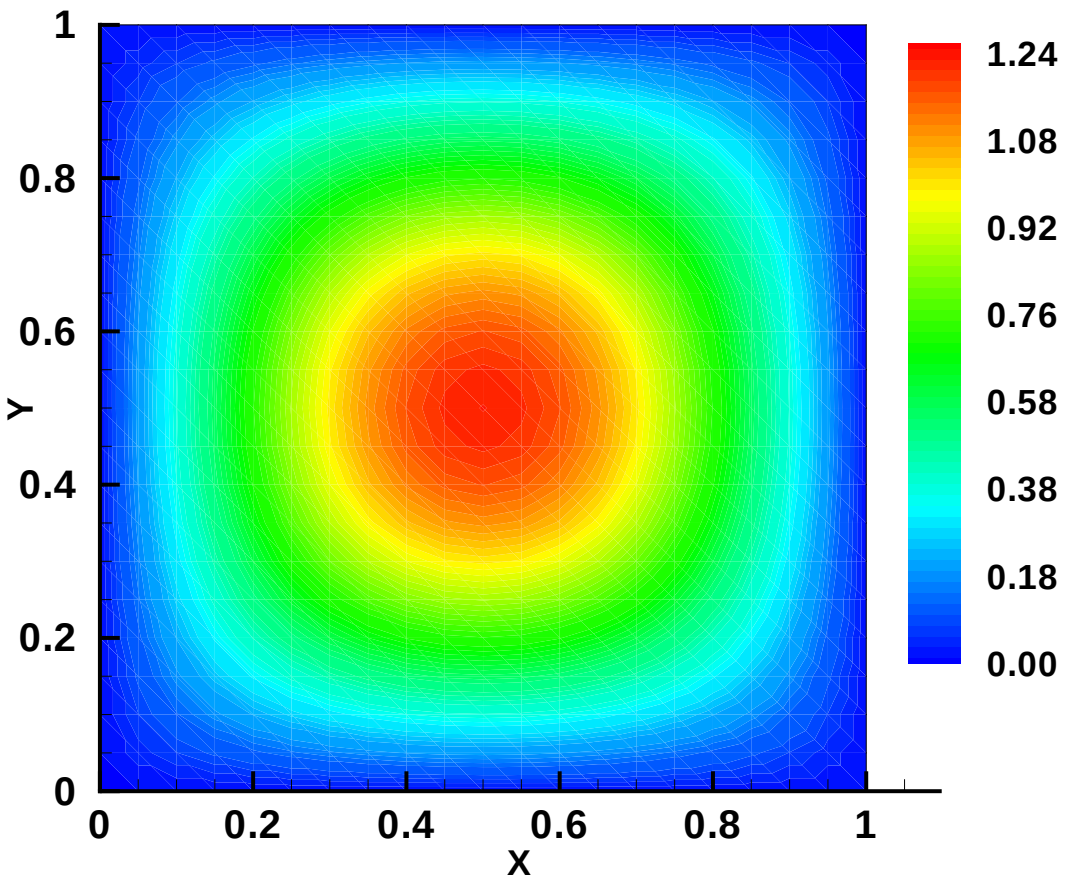

(b) $p$ (analytical)

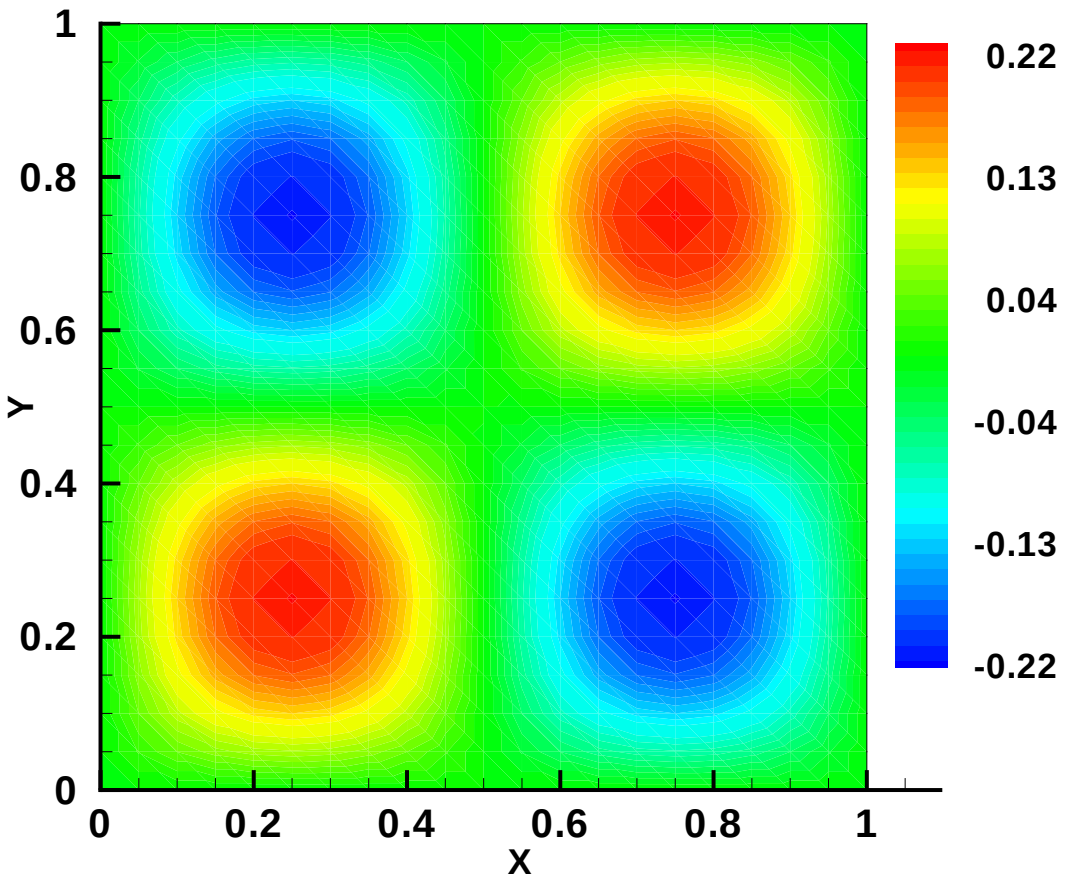

(c) $v_x$ (numerical)

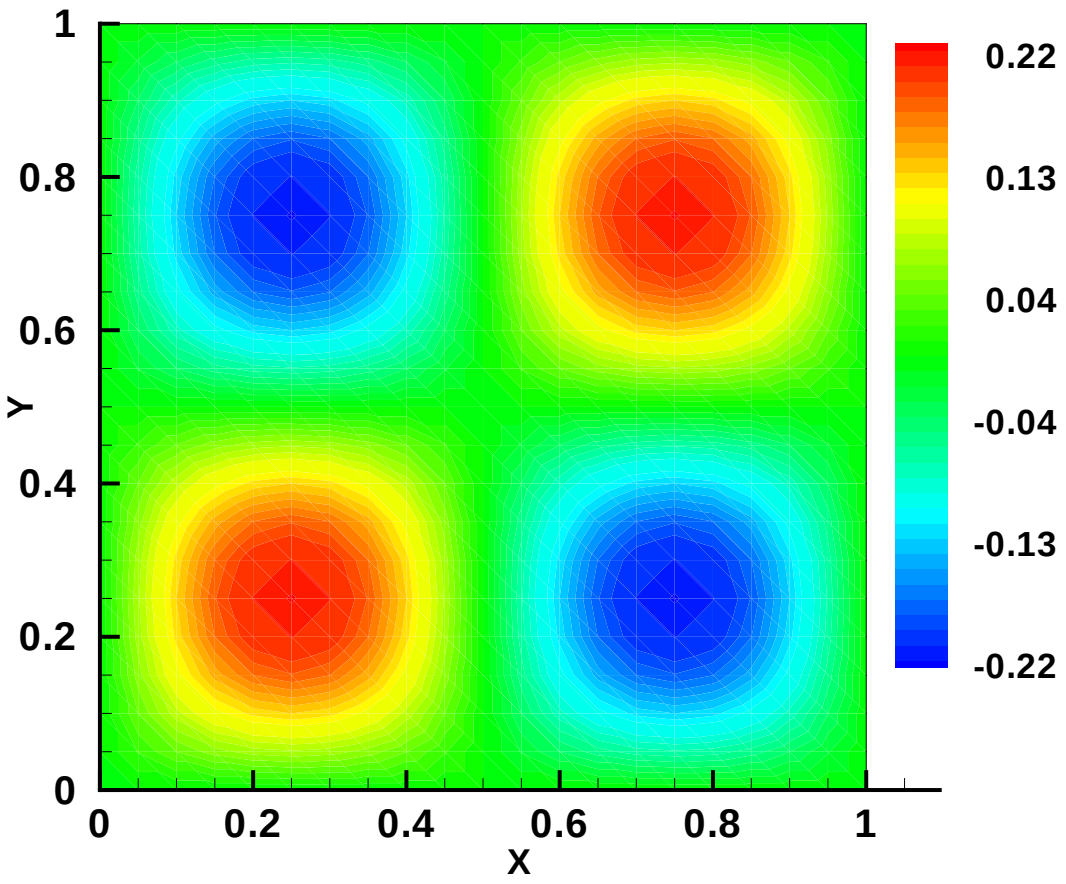

(d) $v_x$ (analytical)

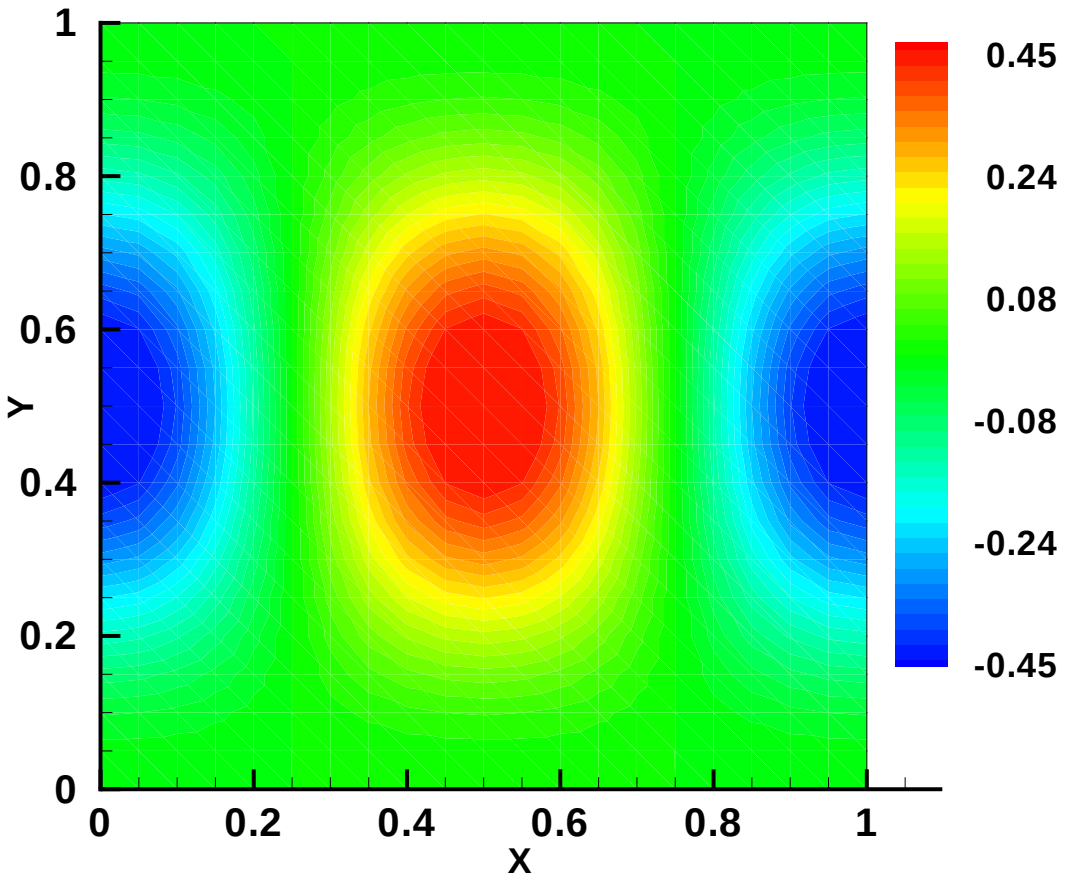

(e) $v_y$ (numerical)

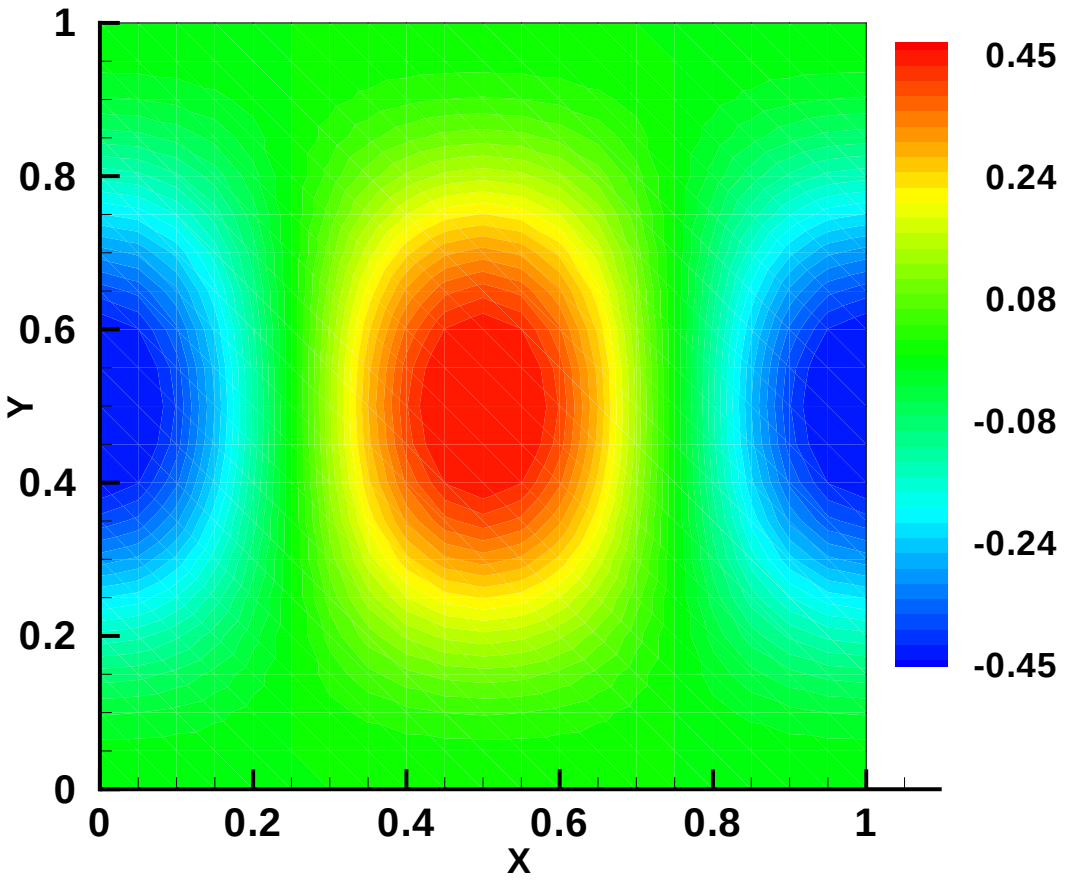

(f) $v_y$ (analytical)

FIGURE 15. Two-dimensional test problem (unsteady Brinkman equation) with *six-node triangular finite elements*: Comparison of analytical and numerical solutions at final time $T = 0.2$. An equally spaced $41 \times 41$ grid was used with $\Delta t = 10^{-3}$.

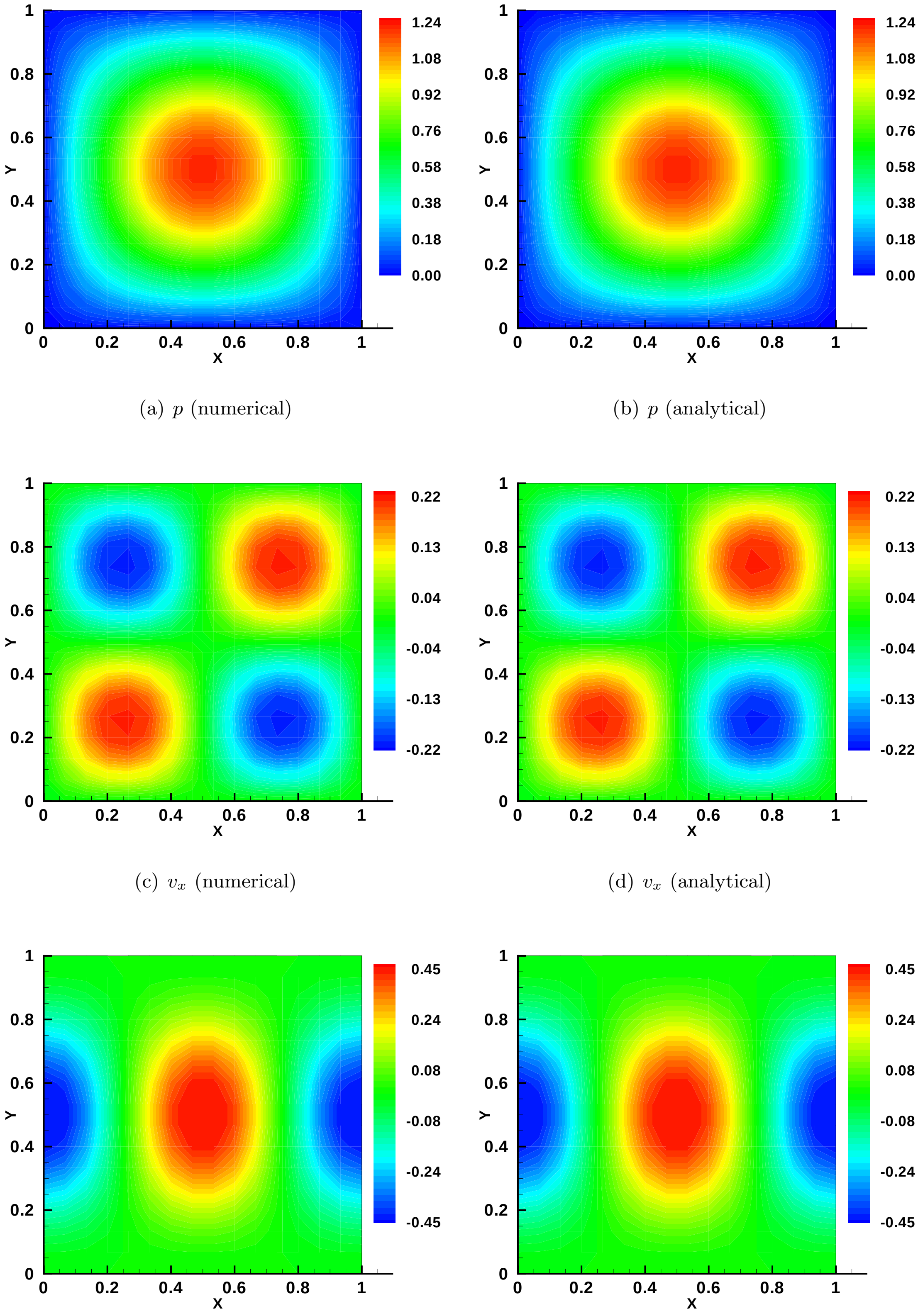

(a) $p$ (numerical) (b) $p$ (analytical)

(c) $v_x$ (numerical) (d) $v_x$ (analytical)

(e) $v_y$ (numerical) (f) $v_y$ (analytical)

FIGURE 16. Two-dimensional test problem (unsteady Brinkman equation) with *nine-node quadrilateral finite elements*: Comparison of analytical and numerical solutions at final time $T = 0.2$. An equally spaced $31 \times 31$ grid was used with $\Delta t = 10^{-3}$.

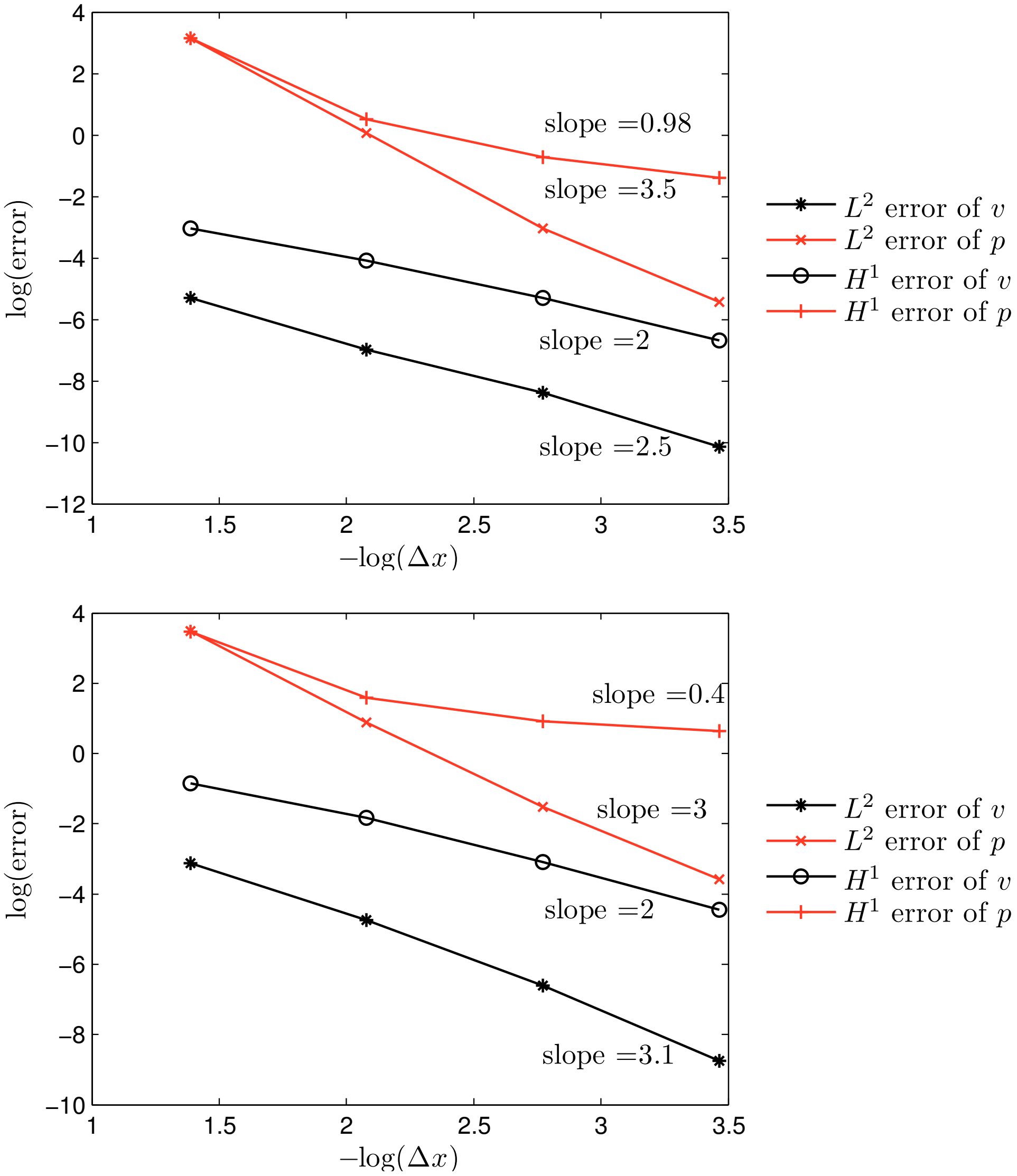

FIGURE 17. Two-dimensional test problem (unsteady Brinkman equation): The figure presents spatial convergence of the proposed formulation using *six-node triangular finite elements*. The time step and the total time are taken to be $\Delta t = 10^{-4}$ and $T = 0.01$ (top) and $\Delta t = 10^{-3}$ and $T = 0.10$ (bottom), respectively.

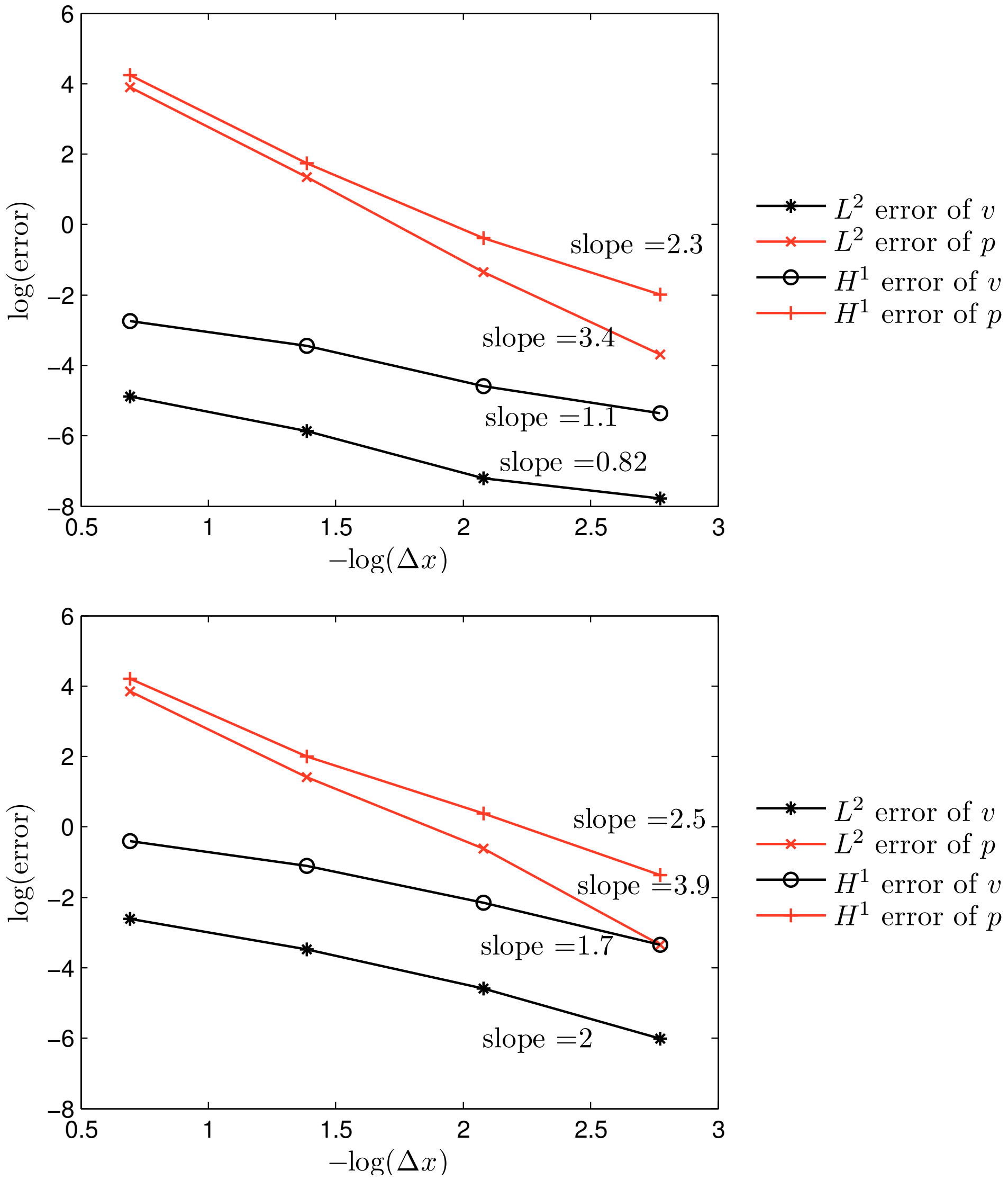

FIGURE 18. Two-dimensional test problem (unsteady Brinkman equation): The figure presents spatial convergence of the proposed formulation using *nine-node quadrilateral finite elements*. The time step and the total time are taken to be $\Delta t = 10^{-4}$ and $T = 0.01$ (top) and $\Delta t = 10^{-3}$ and $T = 0.10$ (bottom), respectively .